\documentclass[smallextended,referee,envcountsect,]{svjour3}
\smartqed
\usepackage{graphicx}
\usepackage{epstopdf}
\usepackage{amsmath}
\usepackage{amssymb}
\usepackage{subfigure}
\usepackage{multicol}
\usepackage{multirow}
\usepackage{float}
\usepackage{makecell}
\usepackage{booktabs}

\journalname{}

\begin{document}

\title{A nonmonotone proximal quasi-Newton method for multiobjective optimization}

%\subtitle{Using  the  LaTex Template}

\author{Xiaoxue Jiang}

\institute{Xiaoxue Jiang \at
            College of Mathematics, Sichuan University \\
               Chengdu 610065, China\\
              xxjiangmath@126.com
}

\date{Received: date / Accepted: date}
%The correct dates will be entered by the editor.

\maketitle

\begin{abstract}
This paper proposes a nonmonotone proximal quasi-Newton algorithm for unconstrained convex multiobjective composite optimization problems.
To design the search direction, we minimize the max-scalarization of the variations of the Hessian approximations and nonsmooth terms.
Subsequently, a nonmonotone line search is used to determine the step size, we allow for the decrease of a convex combination of recent function values.
Under the assumption of strong convexity of the objective function, we prove that the sequence generated by this method converges to a Pareto optimal.
Furthermore, based on the strong convexity, Hessian continuity and Dennis-Mor\'{e} criterion, we use a basic inequality to derive the local superlinear convergence rate of the proposed algorithm. Numerical experiments results demonstrate the feasibility and effectiveness of the proposed algorithm on a set of test problems.
\end{abstract}

%All acknowledgements should be placed in the back of the paper after Conclusions..

\section{Introduction}

%Multiobjective composite optimization problems (MCOPs) are to minimize multiple objective functions simultaneously, each of which consists of the sum of a  smooth and nonsmooth function.
In this paper, we consider the following multiobjective composite optimization problems (MCOPs):
\begin{align}\label{MCOP}
& \min_{x \in \mathbb{R}^{n}}  \ F(x) = (F_{1}(x), F_{2}(x),...,F_{m}(x))^{T}, \tag{MCOP}
\end{align}
where $F_{j}: \mathbb{R}^{n} \rightarrow \mathbb{R}, j = 1, 2, ..., m,$ have the following special separable structure:
$$F_{j}(x) = f_{j}(x) + g_{j}(x),$$
$f_{j}(x)$ is convex and continuously differentiable, the gradient $\nabla f_{j}$ is Lipschitz continuous and $g_{j}(x)$ is convex and continuous but not necessarily differentiable.
Some of the applications of MCOPs include
problems in the areas of image processing \cite{GZY16}, robust optimization \cite{FW14} and machine learning \cite{BMTB20}, etc.

In the past two decades, descent methods have attracted extensive attention in multiobjective optimization community (see, e.g., \cite{AS20,BCNLSS18,FS00,FDS09,GP21,LPP18,P14,WHWLY19} and references therein).
Initially, Fliege et al. \cite{FS00} proposed the multiobjective steepest descent method, which obtains the descent direction by solving an auxiliary and non-parametric scalar subproblem at each iteration, eliminating  the need for predefined parameters.
Since then, the steepest descent method has led to a series of multiobjective optimization algorithms.

Recently, Tanabe et al. \cite{TFY19} proposed a multiobjective proximal gradient algorithm (MPGA) for solving MCOPs.
Subsequently, in \cite{TFY22}, they used a merit function to analyze the convergence rate of MPGA under nonconvex, convex, and strongly convex conditions. These results are consistent with their counterparts in scalar optimization.
The convergence rate analysis in \cite{TFY22} revealed that the MPGA demonstrates slow convergence when encountering ill-conditioned problems.
In order to deal with this problem, Ansary \cite{A22} uses the idea of multiobjective Newton method introduced by Filege et al. \cite{FDS09}  and proposes a Newton-type proximal gradient algorithm (NPGA) for MCOPs, which extends the corresponding method in scalar optimization \cite{LSS14}.
Specifically, in each iteration of NPGA, a subproblem is constructed using the second-order approximation of the smooth term, followed by a Armijo line search strategy to determine the step size.
Under reasonable assumptions, Ansary proved that every accumulation point of the sequence generated by the NGPA converges to the critical point.
However, the convergence rate of the NGPA wasn't obtained in \cite{A22}.
Subsequently, Chen et al. \cite{CJTY23} analyzed the convergence rate of NPGA for MCOPs by using a fundamental inequality, and obtained the quadratic convergence rate of NPGA under the strongly convex assumption.
Likewise, a challenge associated with this method is that NPGA uses precise Hessian information, which is usually not easy to obtain.
In the conclusion of \cite{A22}, Ansary pointed out that NPGA is restricted to convex multiobjective optimization problems and Hessian of each smooth function is used in each iteration which is computationally expensive. And also pointed out that it is a research path to overcome this limitation by using the idea of quasi-Newton method.
In such cases, approximate information can be used to approximate the second order Hessian.

Lately, Peng et al. \cite{PR22} proposed a proximal quasi-Newton algorithm (PQNA) for MCOPs.
The subproblem of PQNA is constructed by using the BFGS quadratic approximation the Hessian of the smooth term and an additional quadratic regularization term.
Then the monotone line search is used to determine the step size.
Under mild assumptions, they obtained that each accumulation point of the sequence generated by PQNA, if exists, is a Pareto stationary point.
However, it's worth mentioning that the superlinear convergence rate of the PQNA wasn't obtained in \cite{PR22}.
It is well known that the single-objective proximal quasi-Newton algorithm leverages the non-expansiveness of the proximal operator to establish superlinear convergence \cite{LSS14}.
But, in the case of multiobjective, when the subproblem is formulated as a proximal operator by involving dual variables, the dual variables are different and variable in each iteration.
Consequently, as in scalar optimization, it is not feasible to prove the superlinear convergence of the multiobjective proximal quasi-Newton algorithm according to the characteristics of the proximal operator.
Meanwhile, in the realm of multiobjective optimization, nonmonotone line search techniques are extensively employed and have been demonstrated to be more effective than monotone algorithms \cite{AS20}.

In this paper, we revisit PQNA for MCOPs.
Taking inspiration from the quasi-Newton method, which approximates second derivative matrices instead of explicitly evaluating them. By incorporating the nonmonotone line search technique, we propose a nonmonotone proximal quasi-Newton algorithm (NPQNA) for MCOPs.
It must be pointed out that the NPQNA of MCOPs differs from the PQNA of MCOPs proposed by Peng et al. in \cite{PR22} in several ways.
Firstly, we use the nonmonotone linear search strategy to find the step size, whereas the latter uses the monotone linear search strategy.
Secondly, if Hessian approximation matrix is positive definite and the subproblem exhibits strong convexity, we do not require the additional quadratic regularization term considered by Peng et al. \cite{PR22} to achieve strong convexity.
Consequently, we omit the regularization term in the construction of the subproblem.
Under standard assumptions, we establish that
the sequence generated by NPQNA converges to a Pareto optimal.
Furthermore, under the assumptions of Hessian continuity, the Dennis-Mor\'{e} criterion and strong convexity, we use a basic inequality to obtain that the local superlinear convergence rate of NPQNA, which was not obtained in Peng et al. \cite{PR22}.
Numerical experiments show that the NPQNA outperform the NPGA \cite{A22} and the PQNA \cite{PR22} in terms of the number of iterations, the number of function evaluations and CPU time for some test problems.

The remainder of this paper is organized as follows. Section 2 gives some preliminary remarks. In Section 3, we recall the nonmonotone line search. Section 4 presents the nonmonotone  proximal quasi-Newton method. In Section 5, we establish the global convergence, and the local superlinear convergence rate of NPQNA. Section 6 presents the results of our numerical experiments.

\section{Preliminaries}

Throughout this paper, the $n$-dimensional Euclidean space $\mathbb{R}^{n}$  is equipped with the inner product $\langle \cdot, \cdot \rangle$ and the induced norm $\| \cdot \|$.
For a positive definite matrix $H$, the notation $\|x\|_{H} = \sqrt{\langle x, Hx \rangle}$ is used to represent the norm induced by $H$ on vector $x$.
For simplicity, we utilize the notation  $\Lambda_{m} = \{1, 2, ..., m\}$ for any $m \in \mathbb{N}$, and define
$$\Delta_{m} := \left\{\lambda : \sum_{j \in \Lambda_{m}} \lambda_{j} = 1, \lambda_{j} \geq 0\right\}.$$
Let $\mathbb{R}^{m}$  be the $m$-dimensional Euclidean space with the partial order ``$\preceq$'' in $\mathbb{R}^{m}$  induced by the Paretian cone $\mathbb{R}_{+}^{m}$, given by $y \preceq z$ (or $z \succeq y$) if and only if
$z - y \in  \mathbb{R}_{+}^{m}$ with its associate relation ``$\prec$'', given by $y \prec z$ (or $z \succ y$) if and only if
$z - y \in \mathbb{R}_{++}^{m}$, where
$$\mathbb{R}_{+}^{m} := \{x \in \mathbb{R}^{m} | x_{i}\geq 0, \ \forall i \in \Lambda_{m}\},$$
and
$$\mathbb{R}_{++}^{m} := \{x \in \mathbb{R}^{m} | x_{i} > 0, \ \forall i \in \Lambda_{m}\}.$$
Since no unique solution which minimizes all objective functions simultaneously exists, we must decide which objective to improve. Hence, the concept of optimality has to be replaced by the concept of Pareto optimality or efficiency.
Given a vector-valued function $F : \mathbb{R}^{n} \rightarrow \mathbb{R}^{m}$, we analyze the NPQNA for finding a Pareto optimum of (MCOP).
A point $x^{*} \in \mathbb{R}^{n}$ is said to be an Pareto optimum of the
$F$ if there does not exist $x \in \mathbb{R}^{n}$ such that $F(x) \leq F (x^{*})$ and $F (x) \neq F (x^{*})$. A
feasible point $x^{*} \in \mathbb{R}^{n}$ is said to be a weakly Pareto optimum of the $F$ if there does not exist $x \in \mathbb{R}^{n}$ such that $F(x) < F(x^{*})$. It is clear that every efficient solution of the $F$ is a weakly Pareto solution, but the converse is not true.
\begin{definition}\label{wendingdian}(\cite{QLGLJ14})
$x^{*} \in \mathbb{R}^{n}$ is a critical point of the (MCOP), if
$$
R(\partial{F_{j}(x^{*};d)}) \cap  (-\mathbb{R}_{++}^{m}) = \emptyset,
$$
where $R(\partial{F_{j}(x^{*};d)})$ denotes the range or image space of the generalized Jacobian of the function $F_{j}$ at $x^{*}$.
\end{definition}
%{\it Proof}  The proof of the first theorem goes here.
%\qed

The classic directional derivative of $F_{j}$ at $x^{*}$ in direction $d$ is defined as
$$F_{j}^{'}(x^{*};d) := \lim\limits_{\alpha \downarrow 0} \frac{F_{j}(x^{*} + \alpha d) - F_{j}(x^{*})}{\alpha}.$$
Definition \ref{wendingdian} shows that
$$
(F_{1}^{'}(x^{*};d), F_{2}^{'}(x^{*};d),..., F_{m}^{'}(x^{*};d)) \notin  -\mathbb{R}_{++}^{m}, \ {\rm for \ all} \ d \in \mathbb{R}^{n}.
$$
That is
\begin{eqnarray}
\max_{j \in \Lambda_{m}} F_{j}^{'}(x^{*};d) \geq 0, \ {\rm for \ all} \ d \in \mathbb{R}^{n}.
\end{eqnarray}

\begin{definition}\label{stronglyconvex}(\cite{A22})
A function $f: \mathbb{R}^{n} \rightarrow \mathbb{R}$ is said to be a $m$-strongly convex if for every $x , y \in \mathbb{R}^{n}$ and $\alpha \in [0, 1]$,
$$f(\alpha x + (1 - \alpha) y) \leq \alpha f(x) + (1 - \alpha) f(y) - \frac{1}{2}m \alpha (1 - \alpha) \|x - y\|^{2}.$$
\end{definition}

\begin{definition}[\cite{DD11}]\label{defnq}
Suppose $h : \mathbb{R}^{n} \rightarrow (-\infty, +\infty]$ be a proper function and $x \in {\rm dom} (h)$.
Then subdifferential of $h$ at $x$ is denoted by $\partial{h(x)}$ and defined as
\begin{eqnarray*}
\partial{h(x)} := \{\xi \in \mathbb{R}^{n} | h(y) \geq h(x) + \xi^{\top}(y - x), \ {\rm for \ all} \ y \in \mathbb{R}^{n} \}.
\end{eqnarray*}
If $x \notin  {\rm dom} (h)$ then we define $\partial{h(x)} = \emptyset$.
\end{definition}

Next we review some properties of subdifferentials.
\begin{theorem}[Proposition 2.82, \cite{DD11}]\label{subdifferential compact}
Let $h : \mathbb{R}^{n} \rightarrow (-\infty, \infty]$ be a proper convex function, and assume that $x \in {\rm int}({\rm dom} h)$. Then  $\partial{h(x)}$ is nonempty and bounded. Moreover, if $h$
is continuous at $x \in {\rm dom} (h)$, then $\partial{h(x)}$ is compact.
\end{theorem}

\begin{theorem}[Theorem 2.91, \cite{DD11}]\label{subdifferential add}
Consider two proper convex functions $h_{i} : \mathbb{R}^{n} \rightarrow \mathbb{R},
i = 1, 2$. Suppose that ${\rm ri} \ {\rm dom} (h_{1}) \cap {\rm ri} \ {\rm dom} (h_{2}) \neq \emptyset$. Then
\begin{eqnarray*}
\begin{aligned}
\partial({h_{1}(x) + h_{2}(x)}) = \partial{h_{1}(x)} + \partial{h_{2}(x)},
\end{aligned}
\end{eqnarray*}
for every $x \in {\rm dom} (h_{1} + h_{2})$.
\end{theorem}

\begin{theorem}[Theorem 2.96, \cite{DD11}]\label{subdifferential convexhull}
Consider convex functions $h_{j} : \mathbb{R}^{n} \rightarrow \mathbb{R}, j \in \Lambda_{m}$ and let $h(x) = \max\{h_{1}(x), h_{2}(x), . . . , h_{m}(x)\}$. Then
\begin{eqnarray*}
\begin{aligned}
\partial{h(x)} = Co\mathop{\cup}\limits_{j\in I(x)} \partial{h_{j}(x)},
\end{aligned}
\end{eqnarray*}
where $I(x) = \{j \in \Lambda_{m} | h(x) = h_{j}(x)\}$ is the active index set and $Co$ is the convex hull.
\end{theorem}

According to the properties of the subdifferential, we recall the necessary conditions for the critical point of (MCOP).

\begin{lemma}[\cite{A22} \emph{Lemma 1}]\label{exsolution}
Suppose $0 \in Co\mathop{\cup}\limits_{j\in \Lambda_{m}} \partial{F_{j}(x^{*})}$ for some $x^{*} \in \mathbb{R}^{n}$, then $x^{*}$ is a critical point of the (MCOP).
\end{lemma}

\begin{proposition}\label{weakparetooptimal}\cite{FDS09}
The following statements hold:
%Let $F$ is continuously differentiable on $\mathbb{R}^{n}$.
\begin{enumerate}
  \item If $x^{*}$ is locally weakly Pareto optimal, then $x^{*}$ is a critical point for $F$.
  \item If $F$ is  $\mathbb{R}_{+}^{m}$-convex and $x^{*}$ is critical for $F$, then $x^{*}$ is weakly Pareto optimal.
  \item If $F$ is  $\mathbb{R}_{+}^{m}$-strictly convex and $x^{*}$ is critical for $F$, then $x^{*}$ is Pareto optimal.
\end{enumerate}
\end{proposition}

Define
\begin{eqnarray*}
\begin{aligned}
Q_{j}(x,d) &:=  \langle \nabla f_{j}(x), d \rangle + \frac{1}{2}d^{T}B_{j}(x)d + g_{j}(x + d) - g_{j}(x), j \in \Lambda_{m}.
\end{aligned}
\end{eqnarray*}
The matrix $B_{j}(x)$ aims to approximate the Hessian matrix of the function $f_{j}$ at iteration $x$.
Define
$$Q(x,d):= \mathop{\max}\limits_{j \in \Lambda_{m}}Q_{j}(x,d).$$
For any fixed $x \in \mathbb{R}^{n}$, $Q_{j}$ is continuous for $d$,
and hence $Q$ is continuous for $d$.
%From Lemma \ref{exsolution}, $\partial_{d}Q_{j}(x,d)$ is nonempty and bounded.
Denote $$I(x,d):=\{j \in \Lambda_{m} \ | \ Q(x,d) = Q_{j}(x,d)\}.$$
For $x \in \mathbb{R}^{n}$, we solve the following subproblem to find a suitable descent direction of the (\ref{MCOP})

\begin{eqnarray}\label{P(x)}
\begin{aligned}
\mathop{\min}_{d \in \mathbb{R}^{n}}Q(x,d) \Leftrightarrow \mathop{\min}_{d \in \mathbb{R}^{n}}\mathop{\max}_{j \in \Lambda_{m}}Q_{j}(x,d).
\end{aligned}
\end{eqnarray}
Denote
\begin{eqnarray}\label{Sol}
\begin{aligned}
d(x) &= \mathop{\arg\min}\limits_{d \in \mathbb{R}^{n}}Q(x,d) \\
&= \mathop{\arg\min}\limits_{d \in \mathbb{R}^{n}}\max_{j \in \Lambda_{m}} \left\{\nabla f_{j}(x)^{T} d + \frac{1}{2}d^{T}B_{j}(x)d + g_{j}(x + d) - g_{j}(x) \right\},
\end{aligned}
\end{eqnarray}
and
\begin{eqnarray}\label{Val}
\begin{aligned}
\theta(x) &= Q(x,d(x)) \\
&= \max_{j \in \Lambda_{m}} \left\{\nabla f_{j}(x)^{T} d(x) + \frac{1}{2}d(x)^{T}B_{j}(x)d(x) + g_{j}(x + d(x)) - g_{j}(x) \right\}.
\end{aligned}
\end{eqnarray}

\begin{remark}
If $g_{j} = 0$ for all $j \in \Lambda_{m}$, then (\ref{P(x)}) coincides with the subproblem in \cite{AS20}, if $m = 1$ then (\ref{P(x)}) coincides with the subproblem used in Step 3 of Algorithm 1 in \cite{LSS14}, and if $B_{j}(x) = \nabla^{2}f_{j}(x)$ then (\ref{P(x)}) coincides with the subproblem used in Step 2 of Algorithm 1 in \cite{A22}.
\end{remark}

We can see that if the matrix $B_{j}(x) \succeq m I$  for all $j$, then $Q_{j}(x, d)$ is $m$-strongly convex function for $d$. Hence $Q(x, d)$ is $m$-strongly convex function for $d$. This implies that (\ref{P(x)}) has a unique finite minimizer.

Clearly for every $x \in \mathbb{R}^{n}$,
\begin{eqnarray}\label{Val<0}
\begin{aligned}
\theta(x) = Q(x,d(x)) \leq Q(x,0) = 0.
\end{aligned}
\end{eqnarray}
Since $d(x)$ is the solution of (\ref{P(x)}), according to the Lemma \ref{exsolution}, we have
$$0 \in \partial_{d}Q_{j}(x,d(x)).$$
The original problem (\ref{P(x)}) can be transformed into:
\begin{eqnarray}\label{Optim-Model}
\begin{aligned}
& \min \  t \\
& \ {\rm s.t.} \ \nabla f_{j}(x)^{T} d + \frac{1}{2}d^{T}B_{j}(x)d + g_{j}(x + d) - g_{j}(x) - t \leq 0, \ j \in \Lambda_{m}, \\
%& \ \ \ \ \ \  j = 1,...,m \\
& \ \ \ \ \ \ (t,d) \in \mathbb{R} \times \mathbb{R}^{n},
\end{aligned}
\end{eqnarray}
which is a convex quadratic problem, we assume that the $B_{j}(x), j \in \Lambda_{m}$ are positive definite.

The Lagrangian function of the problem (\ref{Optim-Model}) is:
\begin{eqnarray*}
\begin{aligned}
L(t,d,\lambda) = t + \sum\limits_{j \in \Lambda_{m}} \lambda_{j} \left(\nabla f_{j}(x)^{T} d + \frac{1}{2}d^{T}B_{j}(x)d + g_{j}(x + d) - g_{j}(x) - t \right).
\end{aligned}
\end{eqnarray*}
As reported in \cite{GP21}, in the smooth case, there are the well-known Karush-Kuhn-Tucker (KKT)
conditions, which are based on the gradients of the objective functions. In case
the objective functions are nonsmooth, the KKT conditions can be generalized using the concept of subdifferentials. Thus,
%From Karush-Kuhn-Tucker(KKT) condition,
there exist $\lambda \in \Delta_{m}$ and $\xi_{j} \in \partial_{d}{g_{j}(x + d)}$, such that
 \begin{eqnarray}\label{KKT1}
\begin{aligned}
\sum\limits_{j \in \Lambda_{m}} \lambda_{j} \left(\nabla f_{j}(x) + B_{j}(x)d + \xi_{j} \right) = 0,
\end{aligned}
\end{eqnarray}

\begin{eqnarray}\label{KKT4}
\begin{aligned}
\nabla f_{j}(x)^{T} d + \frac{1}{2}d^{T}B_{j}(x)d + g_{j}(x + d) - g_{j}(x) \leq t, \ j \in \Lambda_{m},
\end{aligned}
\end{eqnarray}
\begin{eqnarray}\label{KKT5}
\begin{aligned}
\lambda_{j}\left(\nabla f_{j}(x)^{T} d + \frac{1}{2}d^{T}B_{j}(x)d + g_{j}(x + d) - g_{j}(x) - t \right) = 0, \ j \in \Lambda_{m}.
\end{aligned}
\end{eqnarray}

Thus, if $d(x)$ is  the solution of (\ref{Optim-Model}) and $t = \theta(x)$, then there exists $\lambda \in \triangle_{m}$ such that $(d(x), \theta(x))$ satisfy the condition (\ref{KKT1})-(\ref{KKT5}).

\begin{theorem}\label{The equivalence of solutions}
Suppose
%$f_{j}$ is strongly convex function for all $j$
 $B_{j}(x)$ is a positive definite matrix for all $x \in \mathbb{R}^{n}$
and consider $\theta(x)$ as defined by equality  (\ref{Val}). Then,

(a)For all $x \in \mathbb{R}^{n}$, $\theta(x) \leq 0$.

(b)The following conditions are equivalent:

\ \ (1) The point $x$ is not critical.

\ \ (2) $\theta(x) < 0$.

\ \ (3) $d(x) \neq 0$.

(c)The function $\theta$ is continuous.
\end{theorem}
{\it Proof}
Note that, by equality  (\ref{Val}),
 \begin{eqnarray*}
\begin{aligned}
\theta(x) &= \min\limits_{d \in \mathbb{R}^{n}}\max_{j \in \Lambda_{m}} \nabla f_{j}(x)^{T} d + \frac{1}{2}d^{T}B_{j}(x)d + g_{j}(x + d)) - g_{j}(x) \\
&\leq \max_{j \in \Lambda_{m}} \nabla f_{j}(x)^{T} 0 + \frac{1}{2}0^{T}B_{j}(x)0 + g_{j}(x) - g_{j}(x) \\
&= 0.
\end{aligned}
\end{eqnarray*}
So, item (a) holds.

Let us now prove the equivalences of item (b).

First, assume that (1) holds,  i.e.
$R(\partial F_{j}(x)) \bigcap (-\mathbb{R}_{++}^{m}) \neq \emptyset,$ which in turn means that there exists $\widetilde{d} \in \mathbb{R}^{n}$ such that $\partial F_{j}(x)^{T} \widetilde{d} < 0,$
%that is $(\nabla f_{j}(x) + \xi_{j})^{T} \widetilde{d} < 0, \xi_{j} \in \partial g_{j}(x).$
Thus, using equality  (\ref{Val}), for all $t > 0$, we have
\begin{eqnarray*}
\begin{aligned}
\theta(x) &\leq \max_{j \in \Lambda_{m}} \nabla f_{j}(x)^{T} t  \widetilde{d} + \frac{1}{2} t^{2}\widetilde{d}^{T}B_{j}(x)\widetilde{d} + g_{j}(x + t \widetilde{d}) - g_{j}(x) \\
&= t \max_{j \in \Lambda_{m}} \nabla f_{j}(x)^{T} \widetilde{d} + \frac{t}{2}\widetilde{d}^{T}B_{j}(x)\widetilde{d} + \frac{g_{j}(x + t \widetilde{d}) - g_{j}(x)}{t} \\
&\leq t \max_{j \in \Lambda_{m}} \nabla f_{j}(x)^{T} \widetilde{d} + \frac{t}{2}\widetilde{d}^{T}B_{j}(x)\widetilde{d} + \xi_{j}^{T} \widetilde{d}, \ (\xi_{j} \in \partial g_{j}(x + t\widetilde{d})) \\
&= t \max_{j \in \Lambda_{m}} \nabla f_{j}(x)^{T} \widetilde{d} - \nabla f_{j}(x + t\widetilde{d})^{T} \widetilde{d} + \nabla f_{j}(x + t\widetilde{d})^{T} \widetilde{d} + \frac{t}{2}\widetilde{d}^{T}B_{j}(x)\widetilde{d} + \xi_{j}^{T} \widetilde{d}, \\
&\leq t \max_{j \in \Lambda_{m}} L(-t\widetilde{d})\widetilde{d} + \nabla f_{j}(x + t\widetilde{d})^{T} \widetilde{d} + \frac{t}{2}\widetilde{d}^{T}B_{j}(x)\widetilde{d} + \xi_{j}^{T} \widetilde{d}, \\
&= t \max_{j \in \Lambda_{m}} - t L\widetilde{d}^{2} + (\nabla f_{j}(x + t\widetilde{d}) + \xi_{j})^{T} \widetilde{d} + \frac{t}{2}\widetilde{d}^{T}B_{j}(x)\widetilde{d}.
\end{aligned}
\end{eqnarray*}
Therefore, for $t > 0$ sufficiently small, the third term on the right hand side of the above inequality can be viewed as an infinitesimal of $t$, and $\nabla f_{j}$ is Lipschitz continuous, the right hand side of the above inequality is negative and (2) holds.

To prove that (2) implies (3), recall that $\theta(x)$ is the optimal valve of problem (\ref{P(x)}) and so, being negative, the solution to this problem cannot be $d(x) = 0.$

Finally, let us see that (3) implies (1).

For the purpose of contradiction that, let $x$ be a critical point. Since $B_{j}(x)$ is a positive definite matrix for all $x \in \mathbb{R}^{n}$, from $\theta(x) \leq 0,$
\begin{eqnarray}\label{fconvex}
\begin{aligned}
\nabla f_{j}(x)^{T}d(x) + g_{j}(x + d(x)) - g_{j}(x) \leq -\frac{1}{2} d(x)^{T}B_{j}(x)d(x) < 0.
\end{aligned}
\end{eqnarray}
Since $g_{j}$ is convex, for any $\alpha \in (0,1)$
\begin{eqnarray}\label{gconvex}
\begin{aligned}
g_{j}(x + \alpha d(x)) - g_{j}(x) &= g_{j}\left(\alpha (x + d(x)) + (1 - \alpha) x \right) - g_{j}(x) \\
&\leq \alpha g_{j}(x + d(x)) + (1-\alpha)g_{j}(x) - g_{j}(x) \\
&= \alpha (g_{j}(x + d(x)) - g_{j}(x)).\\
\end{aligned}
\end{eqnarray}
From inequality(\ref{fconvex}) and (\ref{gconvex}),
\begin{eqnarray*}
\begin{aligned}
\alpha \nabla f_{j}(x)^{T}d(x) + g_{j}(x + \alpha d(x)) - g_{j}(x)
&\leq \alpha \nabla f_{j}(x)^{T}d(x) + \alpha \left(g_{j}(x + d(x)) - g_{j}(x)\right) \\
&= \alpha \left(\nabla f_{j}(x)^{T}d(x) + g_{j}(x + d(x)) - g_{j}(x)\right) \\
&< 0.
\end{aligned}
\end{eqnarray*}
This implies
\begin{eqnarray*}
\begin{aligned}
\frac{1}{\alpha} \left(\alpha\nabla f_{j}(x)^{T}d(x) + g_{j}(x + \alpha d(x)) - g_{j}(x)\right) < 0.
\end{aligned}
\end{eqnarray*}
Taking the limit $\alpha \rightarrow 0^{+}$ in the above inequality we have
$F^{'}_{j}(x,d(x)) < 0$ for all $j \in \Lambda_{m}.$ This contradicts with $x$ is a critical point. Therefore, $x$ is not a critical point.

We now prove item (c).

It is easy to see that the function $\max\limits_{j \in \Lambda_{m}} \nabla f_{j}(x)^{T} d(x) + \frac{1}{2}d(x)^{T}B_{j}(x)d(x) + g_{j}(x + d(x)) - g_{j}(x)$ is continuous with respect to $x$ and $d(x)$. Therefore, the optimal value function $\theta(x)$ is also continuous from [\cite{B63}, Maximum Theorem]. Moreover, since the optimal set mapping $d(x)$ is unique, $d(x)$ is continuous from [\cite{H73}, Corollary 8.1]
\qed

\begin{lemma}\label{corollary1}
Suppose $B_{j}(x) \succeq m I$,
%Let $B_{j}(x) \succ 0$,
a point $x \in \mathbb{R}^{n}$ is critical for $f$ if and only if $d(x) = 0$, or equivalently, if and only if $\theta(x) = 0.$
\end{lemma}
{\it Proof}
The proof of sufficiency is the previous Theorem \ref{The equivalence of solutions}.

Conversely suppose $d(x) = 0$. Then from condition (\ref{KKT1}) - (\ref{KKT5}), there exists $\lambda \in \Delta_{m}$ such that %$\sum\limits_{j=1}^m{\lambda_{j}} = 1$ and
$\sum\limits_{j \in \Lambda_{m}} (\nabla f_{j}(x) + \xi_{j}) = 0.$ Where $\xi_{j} \in \partial g_{j}(x)$ for $j \in \Lambda_{m}.$ This implies
\begin{eqnarray*}
\begin{aligned}
0 \in Co\mathop{\cup}\limits_{j\in \Lambda_{m}} \partial{F_{j}(x)}.
\end{aligned}
\end{eqnarray*}
Hence, $x$ is a critical point of the (\ref{MCOP}).
\qed

We develop a nonmonotone line search technique to find a suitable step length that allows the objective function values to increase in some iterations. We introduce a classical average-type nonmonotone line search.

\section{Nonmonotone Line Search}

In nonmonotone line search, some growth in the function value is permitted.
The well-known nonmonotone line search is proposed by Zhang and Hager \cite{ZH04}, which takes the average value of the continuous function. We introduce an average-type  nonmonotone line search extended to the multiobjective case in the \cite{MFY22}.
\begin{eqnarray}\label{Fdescend}
\begin{aligned}
F_{j}(x^{k} + \alpha_{k} d^{k}) \leq C_{j}^{k} + \tau \alpha_{k} \theta(x^{k}), \ {\rm for \ all} \ j,
\end{aligned}
\end{eqnarray}
with $\tau \in (0, 1)$ and $C_{j}^{k} \geq F_{j}(x^{k})$,
where
 \begin{eqnarray}\label{averagelinesearch0}
\begin{aligned}
q_{k + 1} = \eta q_{k} + 1,
\end{aligned}
\end{eqnarray}
\begin{eqnarray}\label{averagelinesearch}
\begin{aligned}
C_{j}^{k + 1} = \frac{\eta q_{k}}{q_{k + 1}}C_{j}^{k} +\frac{1}{q_{k + 1}}F_{j}(x^{k+1}).
\end{aligned}
\end{eqnarray}

We prove that the average-type nonmonotone line search for  proximal quasi-Newton method of MCOPs is well defined.

\begin{theorem}\label{linesearchwelldefined}
Suppose $B_{j}(x) \succeq m I$ for every $j$  and $d(x)$ is the solution of (\ref{P(x)}). Then
$$\theta(x) \leq -\frac{m}{2}\| d(x)\|^{2}.$$
Further, if $x$ is noncritical, the inequality (\ref{Fdescend}) holds for every $\alpha > 0$ sufficiently small.
\end{theorem}
{\it Proof}
Suppose $d(x)$ is the solution of (\ref{P(x)}) and $\theta(x) = Q(x, d(x))$.

Then there exists $\lambda \in \Delta_{m}$ such that $(d(x), \theta(x), \lambda)$ satisfies conditions (\ref{KKT1}) - (\ref{KKT5}). Since $g_{j}$ is convex and $\xi_{j} \in \partial g_{j}(x + d(x)),$ we have
\begin{eqnarray}\label{gxiconvex}
\begin{aligned}
g_{j}(x + d(x)) - g_{j}(x) \leq \xi_{j}^{T} d(x).
\end{aligned}
\end{eqnarray}
Then, multiplying both sides of (\ref{KKT1}) by $d(x)$, we obtain
$$\sum\limits_{j\in \Lambda_{m}} \lambda_{j} \left(\nabla f_{j}(x)^{T} d(x) + d(x)^{T}B_{j}(x)d(x) + \xi_{j}^{T}d(x)\right) = 0.$$
Hence from inequality (\ref{gxiconvex}), we get
\begin{eqnarray*}
\begin{aligned}
&\sum\limits_{j\in \Lambda_{m}} \lambda_{j} \left(\nabla f_{j}(x)^{T} d(x) + d(x)^{T}B_{j}(x)d(x) + g_{j}(x + d(x)) - g_{j}(x)\right)  \\
&\leq  \sum\limits_{j\in \Lambda_{m}} \lambda_{j} \left(\nabla f_{j}(x)^{T} d(x) + d(x)^{T}B_{j}(x)d(x) + \xi_{j}^{T}d(x)\right) \\
&= 0. \\
\end{aligned}
\end{eqnarray*}
Take sum over $j \in \Lambda_{m}$ in  (\ref{KKT5})
and using $\sum\limits_{j\in \Lambda_{m}}{\lambda_{j}} = 1$, we have
 \begin{eqnarray*}
\begin{aligned}
  &\sum\limits_{j\in \Lambda_{m}} \lambda_{j} \left(\nabla f_{j}(x)^{T}d(x) + d(x)^{T}B_{j}(x)d(x) + g_{j}(x + d(x))- g_{j}(x)\right) \\
  &=  \frac{1}{2}\sum\limits_{j\in \Lambda_{m}}{\lambda_{j}}d(x)^{T}B_{j}(x)d(x) +  \theta(x) \\
  &\leq 0.
\end{aligned}
\end{eqnarray*}
Using the above two formulas can be obtained
\begin{eqnarray}\label{theta1}
\begin{aligned}
\theta(x) \leq -\frac{1}{2}\sum\limits_{j\in \Lambda_{m}}{\lambda_{j}}d(x)^{T}B_{j}(x)d(x).
\end{aligned}
\end{eqnarray}
Since $B_{j}(x) \succeq m I$ for every $j$, $d(x)^{T}B_{j}(x)d(x) \geq m \| d(x)\|^{2}$ holds for every $j$.
Then from inequality (\ref{theta1}) and $\sum\limits_{j\in \Lambda_{m}}{\lambda_{j}} = 1$, we have
\begin{eqnarray*}
\begin{aligned}
\theta(x) \leq -\frac{m}{2} \left\| d(x) \right\|^{2}.
\end{aligned}
\end{eqnarray*}

Suppose $x$ is noncritical, then from Theorem \ref{The equivalence of solutions}, we get $d(x) \neq 0$. So,
$$\theta(x) \leq -\frac{m}{2}\| d(x)\|^{2} < 0.$$
Since $g_{j}$ is convex, for any $\alpha \in [0,1]$, we have
\begin{eqnarray*}
\begin{aligned}
F_{j}(x + \alpha d(x)) - F_{j}(x)
&= f_{j}(x + \alpha d(x)) - f_{j}(x) + g_{j}(x + \alpha d(x)) - g _{j}(x) \\
&\leq \alpha \nabla f_{j}(x)^{T}d(x) + (1 - \alpha)g _{j}(x) + \alpha g_{j}(x + d(x)) \\
&\ \ \ - g _{j}(x) + o(\alpha^{2}) \\
&= \alpha \left(\nabla f_{j}(x)^{T}d(x) + g_{j}(x + d(x)) - g _{j}(x)\right) + o(\alpha^{2}) \\
&\leq \alpha \theta(x) + o(\alpha^{2})  \ \ \ \ (\frac{1}{2}d(x)^{T}B_{j}(x)d(x) > 0).
\end{aligned}
\end{eqnarray*}
Then for each $j \in \Lambda_{m}$, subtract $\tau\alpha\theta(x)$ on both sides of the about inequality, we obtain
 \begin{eqnarray}\label{linesearchinequality}
\begin{aligned}
F_{j}(x + \alpha d(x))- C_{j}^{k} - \tau\alpha\theta(x) &< F_{j}(x + \alpha d(x)) - F_{j}(x) - \tau\alpha\theta(x) \\
&< \alpha(1 - \tau)\theta(x) + o(\alpha^{2}) \\
&< 0.
\end{aligned}
\end{eqnarray}
Since $\tau \in (0, 1)$ and $\theta(x) < 0$, the right hand side term in inequality (\ref{linesearchinequality}) becomes nonpositive for every $\alpha > 0$ sufficiently small. This implies that the average-type Armijo condition holds for every $\alpha > 0$ sufficiently small, and the proof is complete.
\qed

\section{Nonmonotone proximal quasi-Newton method for MCOPs}

The framework for a complete NPQNA for MCOPs is described as follows:

\vspace{0.5em}
\resizebox{\textwidth}{35mm}{
\begin{tabular}{lllllll}
\hline
{\bf Algorithm 1} NPQNA for MCOPs \\
\hline
   {\bf Step 1} Choose initinal approximation $x^{0} \in \mathbb{R}^{n}$, scalars,
  $\rho, \eta, \tau \in (0, 1),$  \\
  and $\epsilon > 0,$ $\mu > 0$,  $q_{0} = 1$, $C_{j}^{0}  = F_{j}(x^{0}),$ $B_{j}(x^{0}) = I, j \in \Lambda_{m}$,\\
    $0 < \eta_{min} \leq \eta_{max} < 1,$ set $k := 0.$ \\
   {\bf Step 2} Solve the subproblem (\ref{P(x)}) to find $d^{k}$ and $\theta(x^{k}).$ \\
   Where, for $k \geq 1$, \\
  $B_{j}(x^{k})=
  \left\{\begin{array}{lc}
  B_{j}(x^{k - 1}) - \frac{B_{j}(x^{k - 1})S_{k - 1}S_{k - 1}^{T}B_{j}(x^{k - 1})}{S_{k - 1}^{T}B_{j}(x^{k - 1})S_{k - 1}} + \frac{y_{j}^{k - 1}(y_{j}^{k - 1})^{T}}{S_{k - 1}^{T}y_{j}^{k - 1}}, & S_{k - 1}^{T}y_{j}^{k - 1} > 0, \\
  B_{j}(x^{k - 1}), & otherwise.\\
  \end{array}
   \right.$\\
   With $S_{k - 1} = x^{k} - x^{k - 1}$ and $y_{j}^{k - 1} = \nabla f_{j}(x^{k}) - \nabla f_{j}(x^{k - 1}).$ \\
  {\bf Step 3} If $|\theta(x^{k})| < \epsilon$, then stop. Else, go to step 4. \\
   {\bf Step 4} Take $\alpha_{k} = \mu\rho^{h_{k}}$, where $h_{k}$ is the smallest nonnegative integer such \\
   that the average-type condition (\ref{Fdescend}) holds. \\
  Where choose $\eta_{min} \leq  \eta \leq \eta_{max}$, the update of $q_{k}$ and $C_{j}^{k}$ is in (\ref{averagelinesearch0}) and (\ref{averagelinesearch}).\\
   %$q_{k + 1} = \eta q_{k} + 1$ \\
%   $C_{j}^{k + 1} = \frac{\eta q_{k}}{q_{k + 1}}C_{j}^{k} +\frac{1}{q_{k + 1}}F_{j}(x^{k+1})$ \\
   {\bf Step 5} Update $x^{k + 1} = x^{k} + \alpha_{k}d^{k}$. \\
   {\bf Step 6} Set $k = k + 1$ and return to Step 2. \\
\hline
\end{tabular}}
\vspace{0.5em}

\section{Convergence Analysis}

Before proving the convergence of NPQNA for MCOPs, we give the following related  technical lemma.
\begin{lemma}\label{max-average}\cite{MFY22}
Let $\{x^{k}\}$ be the sequence generated by Algorithm 1, and $F_{j}$ is bounded from below. Then, $\{C_{j}^{k}\}$ is nonincreasing and admits a limit when $k \rightarrow \infty.$
\end{lemma}

\begin{lemma}\label{Gamma1}\cite{MFY22}
For a sequence of iterates $\{x^{k}\}$ and search directions $\{d^{k}\}$. Then there exist positive constants $\Gamma_{1} \in (0, 1)$ such that
 \begin{eqnarray}\label{gamma1}
\begin{aligned}
\max_{j \in \Lambda_{m}} \{ \nabla f_{j}(x^{k})^{T} d^{k} + g_{j}(x^{k} + d^{k}) - g_{j}(x^{k})\} \leq -\Gamma_{1}\vert \theta(x^{k}) \vert.
\end{aligned}
\end{eqnarray}
\end{lemma}

\begin{lemma}\label{average}
Let $\{x^{k}\}$ be the sequence generated by Algorithm 1. Assume that $\eta < 1$,  Then, we have
$${\lim_{k \to \infty}} \alpha_{k}|\theta(x^{k})| = 0.$$
\end{lemma}
{\it proof}
Recall that from Lemma \ref{max-average}. $\{C_{j}^{k}\}_{k}$ admits a limit for $k \to \infty$. From the definition of $C_{j}^{k + 1}$ in (\ref{averagelinesearch}), we get
\begin{eqnarray*}
\begin{aligned}
C_{j}^{k + 1} &= \frac{\eta q_{k}}{q_{k + 1}}C_{j}^{k} +\frac{1}{q_{k + 1}}F_{j}(x^{k + 1}) \\
&\leq \frac{\eta q_{k}}{q_{k + 1}}C_{j}^{k} +\frac{1}{q_{k + 1}}(C_{j}^{k} + \tau \alpha_{k} \theta(x^{k})) \\
&= (\frac{\eta q_{k}}{q_{k + 1}} + \frac{1}{{q_{k + 1}}})C_{j}^{k} + \frac{1}{{q_{k + 1}}}  \tau \alpha_{k} \theta(x^{k}) \\
&= C_{j}^{k} + \frac{1}{{q_{k + 1}}} \tau \alpha_{k} \theta(x^{k}),
\end{aligned}
\end{eqnarray*}
where the above inequality holds from inequality (\ref{Fdescend}) and the last equality follows from inequality (\ref{averagelinesearch0}). Since $\frac{\tau \alpha_{k}}{{q_{k + 1}}} \geq 0$ and $\theta(x^{k}) \leq 0$, since $\alpha_{k} = \rho^{h_{k}}$, as $k \to \infty$, $\alpha_{k} \to 0$, while $q_{k + 1} \in (0, 1)$. So we have
\begin{eqnarray}\label{alphatheta}
\begin{aligned}
\lim_{k \to \infty}\frac{\alpha_{k}}{q_{k + 1}}\theta(x^{k}) = 0.
\end{aligned}
\end{eqnarray}
Furthermore, from inequality (\ref{gamma1}) in Lemma \ref{Gamma1}, we obtain
\begin{eqnarray}\label{inequalityGamma1}
\begin{aligned}
&\frac{\alpha_{k}}{q_{k + 1}}\left(\nabla f_{j}(x^{k})^{T} d^{k} + g_{j}(x^{k} + d^{k}) - g_{j}(x^{k})\right) \\
%&\leq \frac{\alpha_{k}}{q_{k + 1}}\theta(x_{k}) \ \  (\frac{1}{2}d(x^{k})^{T}B_{j}(x^{k})d(x^{k}) > 0)\\
&\leq - \frac{\alpha_{k}}{q_{k + 1}}\Gamma_{1}|\theta(x^{k})| \\
&\leq 0.
\end{aligned}
\end{eqnarray}
Now, observe that inequality (\ref{averagelinesearch0}) yields
$$q_{k + 1} = 1 + \sum\limits_{j = 0}^{k}\eta^{j + 1} \leq \sum\limits_{j = 0}^{+\infty}\eta^{j} = \frac{1}{1 - \eta}.$$
Since $0 < \eta < 1$, we have that $\{\frac{1}{q_{k + 1}}\}$ is bounded from below. So, we can take the limit in inequality (\ref{inequalityGamma1}) and use inequality (\ref{alphatheta}) to show that
\begin{eqnarray*}
\begin{aligned}
\lim_{k \to \infty}\alpha_{k}|\theta(x^{k})| = 0,
\end{aligned}
\end{eqnarray*}
and the conclusion follows.
\qed

We will establish the convergence of the NPQNA for MCOPs.

\begin{theorem}\label{convergence}
Suppose that $\{x^{k}\}$ is a sequence generated by Algorithm 1, $B_{j}(x)\succeq m I$, $f_{j}(x)$ is strongly convex with module $m > 0$ and bounded from below,  $\nabla f_{j}$ is Lipschitz continuous. Then every accumulation point of $\{x^{k}\}$ is a Pareto optimal of the (MCOP).
\end{theorem}
{\it proof}
Since $\nabla f_{j}$ is Lipschitz continuous for every $j$ with Lipschitz constant $L_{j}$, where $L = \max\limits_{j \in \Lambda_{m}} L_{j}$, from the so-called descent Lemma \cite[Proposition A.24]{B99} for any $\alpha$ we have  %$\alpha = \frac{\alpha_{k}}{\rho} (\alpha_{k} \leq \rho)$ we have
\begin{eqnarray}\label{Secondupperbound}
\begin{aligned}
f_{j}(x^{k} + \alpha d^{k}) %&\leq f_{j}(x^{k}) + \alpha \nabla f_{j}(x^{k})^{T}d^{k} + \frac{L_{j}}{2}\alpha^{2}\|d^{k}\|^{2} \\
&\leq f_{j}(x^{k}) + \alpha \nabla f_{j}(x^{k})^{T}d^{k} + \frac{L}{2}\alpha^{2}\|d^{k}\|^{2}, \ \forall j \in \Lambda_{m},
\end{aligned}
\end{eqnarray}
also since $g_{j}$ is convex, for any $\alpha \in (0, 1)$, we have
\begin{eqnarray}\label{gjconvex}
\begin{aligned}
g_{j}(x^{k} + \alpha d^{k}) - g_{j}(x^{k})
&= g_{j}\left(\alpha(x^{k} + d^{k}) + (1-\alpha)x^{k}\right) - g_{j}(x^{k}) \\
&\leq \alpha g_{j}(x^{k} + d^{k}) + (1-\alpha)g_{j}(x^{k}) - g_{j}(x^{k}) \\
&= \alpha \left(g_{j}(x^{k} + d^{k}) - g_{j}(x^{k})\right).\\
\end{aligned}
\end{eqnarray}
From inequality (\ref{Secondupperbound}) and inequality (\ref{gjconvex}) for any $\alpha \in (0, 1)$, we get

\begin{eqnarray*}
\begin{aligned}
F_{j}(x^{k} + \alpha d^{k}) - F_{j}(x^{k})
&= f_{j}(x^{k} + \alpha d^{k}) + g_{j}(x^{k} + \alpha d^{k}) - f_{j}(x^{k}) - g_{j}(x^{k}) \\
&\leq \alpha \nabla f_{j}(x^{k})^{T}d^{k} + \frac{L}{2}\alpha^{2}\|d^{k}\|^{2} + \alpha \left(g_{j}(x^{k} + d^{k}) - g_{j}(x^{k})\right) \\
&\leq \alpha \left(\nabla f_{j}(x^{k})^{T}d^{k} + g_{j}(x^{k} + d^{k}) - g_{j}(x^{k})\right) + \frac{L}{2}\alpha^{2}\| d^{k} \|^{2} \\
&\leq \alpha \left(\nabla f_{j}(x^{k})^{T}d^{k} + \frac{1}{2}(d^{k})^{T}B_{k}d^{k} + g_{j}(x^{k} + d^{k}) - g_{j}(x^{k})\right) \\
&\ \ \ + \frac{L}{2}\alpha^{2}\| d^{k} \|^{2} \\
&\leq  \alpha\theta(x^{k}) + \frac{L}{2}\alpha^{2}\| d^{k} \|^{2}, \ j \in \Lambda_{m}, \\
\end{aligned}
\end{eqnarray*}
where the third inequality holds since $B_{j}(x)\succeq m I$. From Step 4 of Algorithm 1, either $\alpha_{k} \geq \rho$ or there exists $h_{k}$ such that $\alpha_{k} = \rho^{h_{k}} < \rho$ satisfies inequality (\ref{Fdescend}).
Since $F_{j}(x^{k}) \leq C_{j}^{k}$, if $\alpha = \frac{\alpha_{k}}{\rho}$ doesn't satisfy inequality  (\ref{Fdescend}), we can get
\begin{eqnarray*}
\begin{aligned}
F_{j}(x^{k} + \frac{\alpha_{k}}{\rho} d^{k})
&> C_{j}^{k} + \tau\frac{\alpha_{k}}{\rho}\theta(x^{k}) \\
&\geq F_{j}(x^{k}) + \tau\frac{\alpha_{k}}{\rho}\theta(x^{k}), \ j \in \Lambda_{m}.\\
\end{aligned}
\end{eqnarray*}
From the above two inequalities, we have
\begin{eqnarray*}
\begin{aligned}
\tau\frac{\alpha_{k}}{\rho}\theta(x^{k}) \leq \frac{\alpha_{k}}{\rho}\theta(x^{k}) + \frac{L}{2}(\frac{\alpha_{k}}{\rho})^{2}\| d^{k} \|^{2}.
\end{aligned}
\end{eqnarray*}
Simplify the above inequality
\begin{eqnarray*}
\begin{aligned}
(\tau - 1)\theta(x^{k}) \leq \frac{L}{2}\frac{\alpha_{k}}{\rho}\| d^{k}\|^{2}.
\end{aligned}
\end{eqnarray*}
Therefore
\begin{eqnarray*}
\begin{aligned}
\alpha_{k} \geq \frac{2\rho (1 - \tau)}{L} \frac{|\theta(x^{k})|}{\| d^{k}\|^{2}}.
\end{aligned}
\end{eqnarray*}
Furthermore, by inequality (\ref{Fdescend}), we have
\begin{eqnarray*}
\begin{aligned}
F_{j}(x^{k + 1}) \leq C_{j}^{k} + \tau\frac{2\rho (1 - \tau)}{L} \frac{|\theta(x^{k})|}{\| d^{k}\|^{2}}\theta(x^{k}), \ j \in \Lambda_{m}.
\end{aligned}
\end{eqnarray*}
So,
\begin{eqnarray*}
\begin{aligned}
F_{j}(x^{k + 1}) \leq C_{j}^{k} - \frac{2\tau\rho (1 - \tau)}{L} \frac{|\theta(x^{k})|^{2}}{\| d^{k}\|^{2}}, \ j \in \Lambda_{m}.
\end{aligned}
\end{eqnarray*}
By Theorem \ref{linesearchwelldefined}, we have
\begin{eqnarray*}
\begin{aligned}
|\theta(x^{k})| \geq \frac{m}{2}\| d(x^{k})\|^{2}.
\end{aligned}
\end{eqnarray*}
Finally, we have
\begin{eqnarray}\label{aboutbeta}
\begin{aligned}
F_{j}(x^{k + 1}) &\leq C_{j}^{k} - \frac{m}{2} \frac{2\tau\rho (1 - \tau)}{L} |\theta(x^{k})| \\
&= C_{j}^{k} - \frac{m \tau\rho (1 - \tau)}{L} |\theta(x^{k})|, \ j \in \Lambda_{m}.\\
\end{aligned}
\end{eqnarray}

Denote $\frac{m \tau\rho (1 - \tau)}{L}$ as $\beta$, combining the relations (\ref{averagelinesearch0}) and (\ref{averagelinesearch}) with inequality (\ref{aboutbeta}), we obtain:

\begin{eqnarray}\label{eq1}
\begin{aligned}
C_{j}^{k + 1} &= \frac{\eta q_{k}C_{j}^{k} + F_{j}(x^{k + 1})}{q_{k + 1}} \\
&\leq \frac{\eta q_{k}C_{j}^{k} + C_{j}^{k} - \beta|\theta(x^{k})|}{q_{k + 1}} \\
&= C_{j}^{k} - \frac{\beta|\theta(x^{k})|}{q_{k + 1}}, \ \forall j \in \Lambda_{m}.
\end{aligned}
\end{eqnarray}
Since $F_{j}, j \in \Lambda_{m}$ are bounded from below and $F_{j}(x^{k}) \leq C_{j}^{k}$, for all $k$, we conclude that $C_{j}^{k}, j \in \Lambda_{m}$ are bounded from below. It follows from \eqref{eq1}  and Lemma \eqref{max-average} that
\begin{eqnarray}\label{thetainfty}
\begin{aligned}
\sum\limits_{k = 0}^{+\infty} \frac{|\theta(x^{k})|}{q_{k + 1}}  \leq C_{j}^{0} - C^{\infty}_{j} < \infty.
\end{aligned}
\end{eqnarray}
Since $f_{j}$ is strongly convex and continuously differentiable, $g_{j}$ is convex and continuous, the level set $\mathcal{L}_{F}(x^{0}) \subset \{x : F_{j}(x) \leq C_{j}^{0}\}$ is compact.
We can get $\{x^{k}\} \subset \mathcal{L}_{F}(x^{0})$ by the fact that $F_{j}(x^{k}) \leq C_{j}^{k} \leq C_{j}^{0}$.
From the compactness of $\mathcal{L}_{F}(x^{0})$, it follows that $\{x^{k}\}$ has an accumulation point.
Now, suppose that $x^{*}$ is an accumulation point of the sequence $\{x^{k}\}$. So, there exists a  subsequence $\{x^{k}\}_{k \in K}$ which converges to $x^{*}$. Then we can prove $\theta(x^{*}) = 0$, by contradiction. Assume $\theta(x^{*}) < 0$, which implies that there are $\epsilon > 0$ and $\delta_{0} > 0$ so that for all $0 < \delta< \delta_{0}$ and for all $k \in K$ such that $|x^{k} - x^{*}| \leq \delta$, we have
\begin{eqnarray*}
\begin{aligned}
|\theta(x^{k})| \geq \epsilon >0.
\end{aligned}
\end{eqnarray*}
This means that
\begin{eqnarray}\label{qkinequality}
\begin{aligned}
\sum\limits_{k = 0}^{+\infty} \frac{|\theta(x^{k})|}{q_{k + 1}} \geq \sum\limits_{k \in \{k \in K | \|x^{k} - x^{*}\| \leq \delta\}} \frac{\epsilon}{q_{k + 1}}.
\end{aligned}
\end{eqnarray}
According to inequality (\ref{averagelinesearch0})
%$q_{k + 1} = \eta q_{k} + 1$ %= 1 + \sum\limits_{i = 0}^{k}\prod\limits_{l = 0}^{i}\eta^{k - l} \leq k + 1$
and $\eta_{max} < 1$, we have

\begin{eqnarray*}
\begin{aligned}
q_{k + 1} &= 1 + \sum\limits_{i = 0}^{k}\prod\limits_{l = 0}^{i}\eta^{k - l} \leq 1 + \sum\limits_{i = 0}^{k}\eta^{i}_{max} \\
&\leq \sum\limits_{i = 0}^{+ \infty}\eta^{i}_{max} = \frac{1}{1 - \eta_{max}}.
\end{aligned}
\end{eqnarray*}
Consequently, following inequality (\ref{qkinequality}), we have
\begin{eqnarray*}
\begin{aligned}
\sum\limits_{k = 0}^{+\infty} \frac{|\theta(x^{k})|}{q_{k + 1}}
&\geq \sum\limits_{k \in \{k \in K | \|x^{k} - x^{*}\| \leq \delta\}} \frac{\epsilon}{q_{k + 1}} \\
&\geq \sum\limits_{k \in \{k \in K | \|x^{k} - x^{*}\| \leq \delta\}} (1 - \eta_{max})\epsilon \\
&= +\infty.
\end{aligned}
\end{eqnarray*}
This contradicts inequality (\ref{thetainfty}). Hence, we have $\theta(x^{*}) = 0$ and it follows from Theorem \ref{The equivalence of solutions} and Lemma \ref{corollary1}  that $x^{*}$ is critical for $F$, since $F$ is $\mathbb{R}^{m}_{+}$-strongly convex, on the basis of Proposition \ref{weakparetooptimal}, so  $x^{*}$ is Pareto optimal for $F$.
\qed

We can get the strong convergence of NPQNA for MCOPs.

\begin{theorem}
Suppose $f_{j}$ is strongly convex with module $m > 0$ for all $j \in \Lambda_{m}$. Let $\{x^{k}\}$ be the sequence generated by Algorithm 1. Then $\{x^{k}\}$ converges to some Pareto optimal $x^{*}$.
\end{theorem}
{\it proof}

According to Theorem \ref{convergence}, we have that $x^{*}$ is a Pareto optimal of $F$.
%and $\lim\limits_{k \rightarrow \infty} F(x^{k}) = F(x^{*})$.
By Lemma \ref{max-average}, we have $\{C^{k}\}$ is nonincreasing and admits a limit when $k \rightarrow \infty.$ Let $\lim\limits_{k \rightarrow \infty} C^{k} = F^{*}$.
According to the definition of $C^{k}$ in (\ref{averagelinesearch}) and the definition of $q_{k}$ in (\ref{averagelinesearch0}) , we have
$$F(x^{k}) = q_{k}C^{k} - \eta q_{k - 1}C^{k - 1},$$
Since the limits of $\{C^{k}\}$ exist, we take the limits on both sides of the above equation and obtain
\begin{eqnarray*}
\begin{aligned}
\lim\limits_{k \rightarrow \infty} F(x^{k}) &= \lim\limits_{k \rightarrow \infty} q_{k}C^{k} - \eta q_{k - 1}C^{k - 1}  \\
&= \frac{1}{1 - \eta} F^{*} - \frac{\eta}{1 - \eta} F^{*} \\
&= F^{*}
\end{aligned}
\end{eqnarray*}
Next, we prove the uniqueness of $x^{*}$ by proof by contradiction.
Suppose instead that there is another accumulation point $x_{1}^{*}$.
Due to the strong convexity of $F$,  the following inequality holds:
$$F(\lambda x^{*} + (1 - \lambda) x_{1}^{*}) \prec \lambda F(x^{*}) + (1 - \lambda) F(x_{1}^{*}) = F^{*},$$
where the equality follows from the convergence of $\{F(x^{k})\}$. However, this contradicts the fact that $x^{*}$ is a Pareto optimal. The uniqueness of accumulation point of $\{x^{k}\}$ implies that $\{x^{k}\}$ converges to $x^{*}$.
\qed

Referring to the proof of the convergence rate of the Newton-type proximal gradient method proposed by chen et al. in \cite{CJTY23}, we establish the local superlinear convergence of the NPQNA for MCOPs.

In the same way, we denote by
$$h_{\lambda}(x):= \sum\limits_{j \in \Lambda_{m}} \lambda_{j}h_{j}(x),$$
$$\nabla h_{\lambda}(x):= \sum\limits_{j \in \Lambda_{m}} \lambda_{j}\nabla h_{j}(x),$$
$$\nabla^{2} h_{\lambda}(x):= \sum\limits_{j \in \Lambda_{m}} \lambda_{j}\nabla^{2} h_{j}(x).$$
According to  Sion's minimax theorem \cite{S58},  there exists $\lambda^{k} \in \Delta_{m}$ such that
$$d^{k} = \mathop{\arg\min}\limits_{d \in \mathbb{R}^{n}} \left\{\nabla f_{\lambda^{k}}(x^{k})^{T} d + \frac{1}{2}d^{T}B_{\lambda^{k}}(x^{k})d + g_{\lambda^{k}}(x^{k} + d) - g_{\lambda^{k}}(x^{k}) \right\},$$
by means of the first-order optimality condition, we have
\begin{eqnarray}\label{first-order optimality condition}
\begin{aligned}
-\nabla f_{\lambda^{k}}(x^{k}) - B_{\lambda^{k}}(x^{k})d^{k} \in \partial g_{\lambda^{k}}(x^{k} + d^{k}).
\end{aligned}
\end{eqnarray}

Similar to the modified fundamental inequality in \cite{CJTY23}, we propose the following inequality.

\begin{proposition}\label{proposition}
Suppose $f_{j}$ is strictly convex for all $j \in \Lambda_{m}$. Let $\{x^{k}\}$ be the sequence generated by Algorithm 1. If $\alpha_{k} = 1$, then there exists $\lambda^{k} \in \Delta_{m}$ such that
\begin{eqnarray}\label{fundamental inequality}
\begin{aligned}
&F_{\lambda^{k}}(x^{k + 1}) - F_{\lambda^{k}}(x) \\
&\leq \frac{1}{2} \|x^{k} - x\|^{2}_{B_{\lambda^{k}}} - \frac{1}{2} \|x^{k + 1} - x^{k}\|^{2}_{B_{\lambda^{k}}} - \frac{1}{2} \|x^{k + 1} - x\|^{2}_{B_{\lambda^{k}}} \\
&\ \ \ + \langle x^{k + 1} - x^{k}, \int_{0}^{1}\int_{0}^{1} \nabla^{2}f_{\lambda^{k}}(x^{k} + st(x^{k + 1} - x^{k}))ds(t(x^{k + 1} - x^{k}))dt \rangle \\
&\ \ \ - \langle x - x^{k}, \int_{0}^{1}\int_{0}^{1} \nabla^{2}f_{\lambda^{k}}(x^{k} + st(x - x^{k}))ds(t(x - x^{k}))dt \rangle.
\end{aligned}
\end{eqnarray}
\end{proposition}
{\it proof}
In the light of the twice continuity of  $f_{j}$, we can deduce that
\begin{eqnarray*}
\begin{aligned}
&F_{j}(x^{k + 1}) - F_{j}(x) \\
&= f_{j}(x^{k + 1}) - f_{j}(x^{k}) - (f_{j}(x) - f_{j}(x^{k})) + g_{j}(x^{k + 1}) - g_{j}(x) \\
&= \langle \nabla f_{j}(x^{k}), x^{k + 1} - x^{k} \rangle + \langle \nabla f_{j}(x^{k}), x^{k} - x \rangle + g_{j}(x^{k + 1}) - g_{j}(x) \\
&\ \ \ + \langle x^{k + 1} - x^{k}, \int_{0}^{1}\int_{0}^{1} \nabla^{2}f_{j}(x^{k} + st(x^{k + 1} - x^{k}))ds(t(x^{k + 1} - x^{k}))dt  \rangle \\
&\ \ \  - \langle x - x^{k}, \int_{0}^{1}\int_{0}^{1} \nabla^{2}f_{j}(x^{k} + st(x - x^{k}))ds(t(x - x^{k}))dt \rangle  \\
&= \langle \nabla f_{j}(x^{k}), x^{k + 1} - x \rangle  + g_{j}(x^{k + 1}) - g_{j}(x) \\
&\ \ \ + \langle x^{k + 1} - x^{k}, \int_{0}^{1}\int_{0}^{1} \nabla^{2}f_{j}(x^{k} + st(x^{k + 1} - x^{k}))ds(t(x^{k + 1} - x^{k}))dt  \rangle  \\
&\ \ \ - \langle x - x^{k}, \int_{0}^{1}\int_{0}^{1} \nabla^{2}f_{j}(x^{k} + st(x - x^{k}))ds(t(x - x^{k}))dt \rangle.
\end{aligned}
\end{eqnarray*}
On the other hand, from (\ref{first-order optimality condition}), we have
$$-\nabla f_{\lambda^{k}}(x^{k}) - B_{\lambda^{k}}(x^{k})d^{k} \in \partial g_{\lambda^{k}}(x^{k} + d^{k}), \lambda^{k} \in \Delta_{m}.$$
This, together with the fact that $\alpha_{k} = 1$, implies
\begin{eqnarray*}
\begin{aligned}
g_{\lambda^{k}}(x^{k + 1}) - g_{\lambda^{k}}(x) &\leq \langle -\nabla f_{\lambda^{k}}(x^{k}) - B_{\lambda^{k}}(x^{k})d^{k}, x^{k + 1} - x \rangle \\
&= \langle -\nabla f_{\lambda^{k}}(x^{k}) - B_{\lambda^{k}}(x^{k})(x^{k + 1} - x^{k}), x^{k + 1} - x \rangle.
\end{aligned}
\end{eqnarray*}
with the help of the last two relations, we have
\begin{eqnarray*}
\begin{aligned}
&F_{\lambda^{k}}(x^{k + 1}) - F_{\lambda^{k}}(x) \\
&\leq \langle B_{\lambda^{k}}(x^{k})(x^{k} - x^{k + 1}) , x^{k + 1} - x \rangle \\
&\ \ \ + \langle x^{k + 1} - x^{k}, \int_{0}^{1}\int_{0}^{1} \nabla^{2}f_{\lambda^{k}}(x^{k} + st(x^{k + 1} - x^{k}))ds(t(x^{k + 1} - x^{k}))dt  \rangle \\
&\ \ \  - \langle x - x^{k}, \int_{0}^{1}\int_{0}^{1} \nabla^{2}f_{\lambda^{k}}(x^{k} + st(x - x^{k}))ds(t(x - x^{k}))dt \rangle  \\
&= \frac{1}{2} \|x^{k} - x\|^{2}_{B_{\lambda^{k}}} - \frac{1}{2} \|x^{k + 1} - x^{k}\|^{2}_{B_{\lambda^{k}}} - \frac{1}{2} \|x^{k + 1} - x\|^{2}_{B_{\lambda^{k}}} \\
&\ \ \ + \langle x^{k + 1} - x^{k}, \int_{0}^{1}\int_{0}^{1} \nabla^{2}f_{\lambda^{k}}(x^{k} + st(x^{k + 1} - x^{k}))ds(t(x^{k + 1} - x^{k}))dt  \rangle \\
&\ \ \  - \langle x - x^{k}, \int_{0}^{1}\int_{0}^{1} \nabla^{2}f_{\lambda^{k}}(x^{k} + st(x - x^{k}))ds(t(x - x^{k}))dt \rangle.  \\
\end{aligned}
\end{eqnarray*}
\qed

Again, assuming that $f_{j}$ is twice continuously differentiable and strongly convex with constant $m$, and  $\nabla^{2} f_{j}$ are continuous, and the following assumption holds.

\newtheorem{assumption}{Assumption}[section]
\begin{assumption}\label{assumption4}\cite{AS20}
For all $\epsilon > 0$ there exist $k^{0} \in \mathbb{N}$ such that for all $k \geq k^{0}$, we have
\begin{eqnarray}\label{assump4}
\begin{aligned}
\frac{\left\| (\nabla^{2} f_{j}(x^{*}) - B_{j}(x^{k}))d^{k} \right\|}{\|d^{k}\|} \leq \epsilon, j \in \Lambda_{m},
\end{aligned}
\end{eqnarray}
where $\{x^{k}\}$ is a sequence generated by Algorithm 1 and $\{B_{j}(x^{k})\}, j \in \Lambda_{m},$ are the sequences obtained by the BFGS updates.
\begin{eqnarray*}
\begin{aligned}
B_{j}(x^{k + 1}) = B_{j}(x^{k}) -  \frac{B_{j}(x^{k}) s^{k}(s^{k})^{T} B_{j}(x^{k})}{(s^{k})^{T} B_{j}(x^{k})s^{k}} + \frac{y_{j}^{k}(y_{j}^{k})^{T}}{(y_{j}^{k})^{T}s^{k}},
\end{aligned}
\end{eqnarray*}
where %$\alpha_{k}$ is a suitable step size and
\begin{eqnarray*}
\begin{aligned}
s^{k} = x^{k + 1} - x^{k}, \ y_{j}^{k} = \nabla f_{j}(x^{k + 1}) - \nabla f_{j}(x^{k}).
\end{aligned}
\end{eqnarray*}
The Assumption \ref{assumption4} is also called Dennis-Mor\'{e} criterion.
\end{assumption}

\begin{lemma}\label{lemma1}
Suppose $f_{j}$ is twice continuously differentiable and strongly convex with constant $m$, $\nabla^{2} f_{j}(x^{k})$ is continuous, the sequence $\{B_{j}(x^{k})\}$ satisfies the Assumption (\ref{assumption4}). Let $\{x^{k}\}$ be the sequence generated by Algorithm 1.  Then, for any $0 < \epsilon \leq \frac{m(1 - \tau)}{3}$, the unit step length satisfies the nonmonotone line search conditions (\ref{Fdescend}) after sufficiently many iterations.
\end{lemma}
{\it proof}

Since $f_{j}$ is twice continuously differentiable, we have
\begin{eqnarray*}
\begin{aligned}
f_{j}(x^{k} + d^{k}) \leq f_{j}(x^{k}) + \nabla f_{j}(x^{k})^{T}d^{k} + \frac{1}{2}(d^{k})^{T}\nabla^{2} f_{j}(x^{k})d^{k} + \frac{\epsilon}{2}\left\|d^{k}\right\|^{2},
\end{aligned}
\end{eqnarray*}
we add $g_{j}(x^{k} + d^{k})$ to both sides to obtain
\begin{eqnarray*}
\begin{aligned}
F_{j}(x^{k} + d^{k}) \leq f_{j}(x^{k}) + \nabla f_{j}(x^{k})^{T}d^{k} + \frac{1}{2}(d^{k})^{T}\nabla^{2} f_{j}(x^{k})d^{k} + \frac{\epsilon}{2}\left\|d^{k}\right\|^{2} + g_{j}(x^{k} + d^{k}),
\end{aligned}
\end{eqnarray*}
we then add and subtract $g_{j}(x^{k})$ from the right-hand side to obtain
\begin{eqnarray}\label{Leminequality}
\begin{aligned}
F_{j}(x^{k} + d^{k}) &\leq f_{j}(x^{k}) + g_{j}(x^{k}) + \nabla f_{j}(x^{k})^{T}d^{k} + g_{j}(x^{k} + d^{k}) - g_{j}(x^{k}) \\
&\ \ \ + \frac{1}{2}(d^{k})^{T}\nabla^{2} f_{j}(x^{k})d^{k} + \frac{\epsilon}{2}\left\|d^{k}\right\|^{2} + \frac{1}{2}(d^{k})^{T}B_{j}(x^{k})d^{k} \\
&\ \ \ - \frac{1}{2}(d^{k})^{T}B_{j}(x^{k})d^{k} \\
&\leq F_{j}(x^{k}) + \theta(x^{k}) + \frac{1}{2}(d^{k})^{T}\left(\nabla^{2} f_{j}(x^{k}) - \nabla^{2} f_{j}(x^{*}))\right)d^{k}  \\
&\ \ \ +  \frac{1}{2}(d^{k})^{T}\left(\nabla^{2} f_{j}(x^{*}) - B_{j}(x^{k})\right)d^{k} +  \frac{\epsilon}{2}\left\|d^{k}\right\|^{2}  \\
&\leq F_{j}(x^{k}) + \theta(x^{k}) +  \frac{3 \epsilon}{2}\left\|d^{k}\right\|^{2} \\
&\leq F_{j}(x^{k}) + \tau \theta(x^{k}) + (1 - \tau)\theta(x^{k}) + \frac{3 \epsilon}{2}\left\|d^{k}\right\|^{2} \\
&\leq F_{j}(x^{k}) + \tau \theta(x^{k}) + (\frac{3 \epsilon}{2} - \frac{m(1 - \tau)}{2}) \left\|d^{k}\right\|^{2} \\
&\leq C_{j}^{k} + \tau \theta(x^{k}),
\end{aligned}
\end{eqnarray}
where the third inequality comes from the Assumption \ref{assumption4} and $\nabla^{2} f_{j}$ is continuous. The second to last inequality is derived from the Theorem \ref{linesearchwelldefined}, and $0 < \epsilon \leq \frac{(1 - \tau)m}{3}$, we have
\begin{eqnarray*}
\begin{aligned}
\frac{3 \epsilon}{2} - \frac{m(1 - \tau)}{2} \leq 0,
\end{aligned}
\end{eqnarray*}
so, the last inequality holds.
therefore, the nonmonotone line search conditions hold for $\alpha_{k} = 1.$
\qed

\begin{theorem}
Suppose $f_{j}$ is strongly convex with module $m > 0$, and its Hessian is continuous
for $j \in \Lambda_{m}$, and $\{B_{j}(x^{k})\}$ satisfies the Assumption \ref{assumption4}  and $B_{j}(x^{k}) \succeq m I $ for some $m > 0$. Let $\{x^{k}\}$ denote a bounded sequence generated by Algorithm 1. Then, for any $0 < \epsilon \leq \frac{m(1 - \tau)}{3}$,
there exists $K_{\epsilon} > 0$ such that
\begin{eqnarray*}
\begin{aligned}
\|x^{k + 1} - x^{*}\| \leq \sqrt{\frac{3 \epsilon (1 + \tau_{k}^{2})}{m}} \|x^{k}- x^{*}\|.
\end{aligned}
\end{eqnarray*}
holds  for all $k \geq K_{\epsilon}$, where $\tau_{k} := \frac{\|x^{k + 1} - x^{k}\|}{\|x^{k} - x^{*}\|} \in [\frac{m - \sqrt{6 m\epsilon - 9 \epsilon^{2}}}{m - 3 \epsilon}, \frac{m + \sqrt{6 m\epsilon - 9 \epsilon^{2}}}{m - 3\epsilon}]$. Furthermore,
the sequence $\{x^{k}\}$ converges superlinearly to $x^{*}$.
\end{theorem}
{\it proof}
Since the assumptions of Lemma \ref{lemma1} are satisfied, unit step lengths satisfy the nonmonotone line search conditions after sufficiently many iterations:
$$\alpha_{k} = 1,$$
$$x^{k+1} = x^{k} + d^{k}.$$
Substituting $x = x^{*}$ into inequality (\ref{fundamental inequality}), we obtain
\begin{eqnarray}\label{fundamental inequality1111}
\begin{aligned}
0 &\leq F_{\lambda^{k}}(x^{k + 1}) - F_{\lambda^{k}}(x^{*}) \\
&\leq \frac{1}{2} \|x^{k} - x^{*}\|^{2}_{B_{\lambda^{k}}(x^{k})} - \frac{1}{2} \|x^{k + 1} - x^{k}\|^{2}_{B_{\lambda^{k}}(x^{k})} - \frac{1}{2} \|x^{k + 1} - x^{*}\|^{2}_{B_{\lambda^{k}}(x^{k})} \\
&\ \ \ + \langle x^{k + 1} - x^{k}, \int_{0}^{1}\int_{0}^{1} \nabla^{2}f_{\lambda^{k}}(x^{k} + st(x^{k + 1} - x^{k}))ds(t(x^{k + 1} - x^{k}))dt \rangle \\
&\ \ \ - \langle x^{*} - x^{k}, \int_{0}^{1}\int_{0}^{1} \nabla^{2}f_{\lambda^{k}}(x^{k} + st(x^{*} - x^{k}))ds(t(x^{*} - x^{k}))dt \rangle.
\end{aligned}
\end{eqnarray}
where the first inequality comes from the fact $F(x^{*}) \preceq F(x^{k})$ for all $k$. On the other hand,
$$\frac{1}{2} \|x - x^{k}\|^{2}_{B_{\lambda^{k}}(x^{k})} = \langle x - x^{k}, \int_{0}^{1}\int_{0}^{1} B_{\lambda^{k}}(x^{k})ds(t(x - x^{k}))dt \rangle,$$
then,  inserting $x = x^{*}$ and $x = x^{k + 1}$ respectively into the above formula, we have
$$\frac{1}{2} \|x^{*} - x^{k}\|^{2}_{B_{\lambda^{k}}(x^{k})} = \langle x^{*} - x^{k}, \int_{0}^{1}\int_{0}^{1} B_{\lambda^{k}}(x^{k})ds(t(x^{*} - x^{k}))dt \rangle,$$
$$\frac{1}{2} \|x^{k + 1} - x^{k}\|^{2}_{B_{\lambda^{k}}(x^{k})} = \langle x^{k + 1} - x^{k}, \int_{0}^{1}\int_{0}^{1} B_{\lambda^{k}}(x^{k})ds(t(x^{k + 1} - x^{k}))dt \rangle.$$
Since the accumulation point of $\{x^{k}\}$ is Pareto optimal, there exists $K_{\epsilon} \geq K_{\epsilon}^{1}$ such that, for all $k \geq K_{\epsilon}$,
$$\|\nabla^{2}f_{\lambda^{k}}(x^{k} + st(x^{k + 1} - x^{k}) - \nabla^{2}f_{\lambda^{k}}(x^{k})\| \leq \epsilon, \forall s, t \in [0, 1],$$
$$\|\nabla^{2}f_{\lambda^{k}}(x^{k} + st(x^{*} - x^{k}) - \nabla^{2}f_{\lambda^{k}}(x^{k})\| \leq \epsilon, \forall s, t \in [0, 1].$$
Substituting the above relation into (\ref{fundamental inequality1111}), we obtain
\begin{eqnarray*}
\begin{aligned}
&\frac{1}{2} \|x^{k + 1} - x^{*}\|^{2}_{B_{\lambda^{k}}(x^{k})}  \\
&\leq \langle x^{k + 1} - x^{k}, \int_{0}^{1}\int_{0}^{1}
(\nabla^{2}f_{\lambda^{k}}(x^{k} + st(x^{k + 1} - x^{k}))  - B_{\lambda^{k}}(x^{k}))ds(t(x^{k + 1} - x^{k}))dt \rangle \\
&\ \ \ - \langle x^{*} - x^{k}, \int_{0}^{1}\int_{0}^{1} (\nabla^{2}f_{\lambda^{k}}(x^{k} + st(x^{*} - x^{k})) - B_{\lambda^{k}}(x^{k}))ds(t(x^{*} - x^{k}))dt \rangle  \\
&\leq \frac{\epsilon}{2}\|x^{k + 1} - x^{k}\|^{2} + \frac{\epsilon}{2}\|x^{*} - x^{k}\|^{2} \\
&\ \ \ + \langle x^{k + 1} - x^{k}, \int_{0}^{1}\int_{0}^{1} (\nabla^{2}f_{\lambda^{k}}(x^{k})- B_{\lambda^{k}}(x^{k}))ds(t(x^{k + 1} - x^{k}))dt \rangle \\
&\ \ \ - \langle x^{*} - x^{k}, \int_{0}^{1}\int_{0}^{1} (\nabla^{2}f_{\lambda^{k}}(x^{k})- B_{\lambda^{k}}(x^{k}))ds(t(x^{*} - x^{k}))dt \rangle \\
&\leq \epsilon \|x^{k + 1} - x^{k}\|^{2} + \epsilon \|x^{*} - x^{k}\|^{2} \\
&\ \ \ + \langle x^{k + 1} - x^{k}, \int_{0}^{1}\int_{0}^{1} (\nabla^{2}f_{\lambda^{k}}(x^{*})- B_{\lambda^{k}}(x^{k}))ds(t(x^{k + 1} - x^{k}))dt \rangle \\
&\ \ \ - \langle x^{*} - x^{k}, \int_{0}^{1}\int_{0}^{1} (\nabla^{2}f_{\lambda^{k}}(x^{*})- B_{\lambda^{k}}(x^{k}))ds(t(x^{*} - x^{k}))dt \rangle \\
&\leq \frac{3 \epsilon}{2} \|x^{k + 1} - x^{k}\|^{2} + \frac{3 \epsilon}{2} \|x^{*} - x^{k}\|^{2}.
\end{aligned}
\end{eqnarray*}
where the last inequality holds because $\nabla^{2}f_{j}$ is continuous and the Assumption (\ref{assumption4}) holds. This together with $B_{j} \succeq m I$, implies
\begin{eqnarray}\label{ine1111}
\begin{aligned}
m \|x^{k + 1} - x^{*}\|^{2} \leq 3 \epsilon \|x^{k + 1} - x^{k}\|^{2} + 3 \epsilon \|x^{*} - x^{k}\|^{2}.
\end{aligned}
\end{eqnarray}
By direct calculation, we have
\begin{eqnarray*}
\begin{aligned}
&3 \epsilon \|x^{k + 1} - x^{k}\|^{2} + 3 \epsilon \|x^{*} - x^{k}\|^{2} \\
&\geq m \|x^{k + 1} - x^{*}\|^{2} \\
&= m \|x^{k + 1}- x^{k} + x^{k} - x^{*}\|^{2} \\
&\geq m \|x^{k + 1}- x^{k}\|^{2} + m \| x^{k} - x^{*}\|^{2} - 2 m \|x^{k + 1}- x^{k}\| \| x^{k} - x^{*}\|.
\end{aligned}
\end{eqnarray*}
Rearranging, we have
\begin{eqnarray*}
\begin{aligned}
2 m \|x^{k + 1}- x^{k}\| \| x^{k} - x^{*}\| \geq (m - 3 \epsilon) \|x^{k + 1} - x^{k}\|^{2} + (m - 3 \epsilon) \| x^{k} - x^{*}\|^{2} ,
\end{aligned}
\end{eqnarray*}
dividing by $\| x^{k} - x^{*}\|^{2}$, we get
\begin{eqnarray*}
\begin{aligned}
2 m \frac{\|x^{k + 1}- x^{k}\|}{\| x^{k} - x^{*}\|}  \geq (m - 3 \epsilon) \frac{\|x^{k + 1} - x^{k}\|^{2}}{\| x^{k} - x^{*}\|^{2}} + (m - 3 \epsilon),
\end{aligned}
\end{eqnarray*}
let $\tau_{k} = \frac{\|x^{k + 1} - x^{k}\|}{\|x^{k} - x^{*}\|}$, so we get
\begin{eqnarray*}
\begin{aligned}
(m - 3 \epsilon) \tau_{k}^{2} - 2 m \tau_{k} + m - 3 \epsilon \leq 0.
\end{aligned}
\end{eqnarray*}
Since $\epsilon \leq \frac{m(1 -  \tau)}{3}$, we deduce that  $\tau_{k} \in [\frac{m - \sqrt{6 m\epsilon - 9 \epsilon^{2}}}{m - 3 \epsilon}, \frac{m + \sqrt{6 m\epsilon - 9 \epsilon^{2}}}{m - 3\epsilon}]$.
Substituting $\tau_{k}$ into relation  (\ref{ine1111}), we derive that
\begin{eqnarray*}
\begin{aligned}
\|x^{k + 1} - x^{*}\| \leq \sqrt{\frac{3 \epsilon (1 + \tau_{k}^{2})}{m}} \|x^{k} - x^{*}\|.
\end{aligned}
\end{eqnarray*}
Furthermore, since $\epsilon$ tends to 0 as $k$ tends to infinity, it follows that
\begin{eqnarray*}
\begin{aligned}
\lim_{k \rightarrow \infty} \tau_{k} \in \lim_{\epsilon \rightarrow 0} [\frac{m - \sqrt{6 m\epsilon - 9 \epsilon^{2}}}{m - 3 \epsilon}, \frac{m + \sqrt{6 m\epsilon - 9 \epsilon^{2}}}{m - 3\epsilon}] = \{1\}.
\end{aligned}
\end{eqnarray*}
We use the relation to get
\begin{eqnarray*}
\begin{aligned}
\lim_{k \rightarrow \infty} \frac{\|x^{k + 1} - x^{*}\|}{ \|x^{k} - x^{*}\|} = 0.
\end{aligned}
\end{eqnarray*}
This concludes the proof.
\qed

\section{Numerical Experiment}

In this section, we presents some numerical experiments in order to illustrate the applicability of our Algorithm NPQNA for MCOPs. Based on it, we compare:
\begin{itemize}
  \item  A nonmonotone proximal quasi-Newton algorithm (NPQNA) for MCOPs.
  \item  A proximal quasi-Newton algorithm (PQNA) for MCOPs in \cite{PR22}.
  \item  A Newton-type proximal gradient algorithm (NPGA) for MCOPs in \cite{A22}.
\end{itemize}

We use python 3.9 to program the Algorithm 1, and run on a computer with CPU Intel Core i7 2.90GHz and 32GB of memory, and the subproblems is solved by the solver 'SLSQP' in python.
$\| d^{k}\| \leq 10^{-6}$  or maximum 300 iterations is considered as stopping criteria.
Solution of a multiobjective optimization problem is not isolated optimum points, but a set of efficient solutions. To generate an approximate set of efficient solutions we have considered multi-start technique. Following steps are executed in this technique. \\
      (1) A set of 100 uniformly distributed random initial points between lower bounds $lb$ and upper bounds $ub$ are considered,  where $lb$, $ub \in \mathbb{R}^{n}$ and $lb < ub$.\\
      (2) Algorithm 1 is executed individually.\\
      (3) We implemented the methods using the nonmonotone line search step size strategy with parameters $\eta = 10^{-4}$.

\textbf{Test Problems}: We have constructed a set of nonlinear 2 or 3 objective optimization problems to illustrate the effectiveness of our proposed NPQNA for MCOPs and found a good pareto front. Some differentiable multiobjective test problems $f_{j}(x)$ are paired with the nonsmooth multiobjective problems listed below. Details are provided in Table 1. The dimension of variables and objective functions are presented in the second columns.

Since the subproblem is a convex nonsmooth quadratic constrained optimization problem, where the objective function is a linear function and the constraint is a quadratic constraint function with a nonsmooth term.
\begin{eqnarray}\label{Optim-Model1}
\begin{aligned}
& \min \ t \\
& s.t. \ \nabla f_{j}(x)^{T} d + \frac{1}{2}d^{T}B_{j}(x)d + g_{j}(x + d) - g_{j}(x) - t \leq 0, \ j \in \Lambda_{m}, \\
& \ \ \ \ \ \ (t,d) \in \mathbb{R} \times \mathbb{R}^{n}.
\end{aligned}
\end{eqnarray}
We refer to the way of dealing with nonsmooth term $g_{j}$ in \cite{AFP23}, when we similar consider the linear programming with box constraints in the nonsmooth part of the original multiobjective optimization problem. For each $j \in J$, we assume that
${\rm dom} (g_{j}) = \{x \in \mathbb{R}^{n} | lb \preceq x \preceq ub\} =: C$, where %$lb, ub \in \mathbb{R}^{n}$ are given in the second columns of Table 1, and
define $g_{j} : C\rightarrow \mathbb{R} $ by
%$$g_{j}(x) = \max_{z \in Z_{j}} <x,z>,$$
\begin{eqnarray}\label{g_j}
\begin{aligned}
& g_{j}(x) = \max_{z \in Z_{j}} \langle x, z \rangle,
\end{aligned}
\end{eqnarray}
where $Z_{j} \in  \mathbb{R}^{n}$ is the uncertainty set. Let $B_{j} \in \mathbb{R}^{n} \times \mathbb{R}$ a nonsingular matrix and $\delta > 0$ be given.
We set
%$$Z_{j} := \{ z \in \mathbb{R}^{n} | -\delta e \preceq B_{j}z \preceq \delta e \},$$
\begin{eqnarray}\label{Z_j}
\begin{aligned}
& Z_{j} := \{ z \in \mathbb{R}^{n} | -\delta e \preceq B_{j}z \preceq \delta e \},
\end{aligned}
\end{eqnarray}
where  $e = (1, . . . , 1)^{T} \in \mathbb{R}^{n}$. Since $Z_{j}$ is a nonempty and compact, $g_{j}(x)$ is well-defined. In our tests, the elements of the matrix $B_{j}$ were randomly chosen between 0 and 1. In turn, given an arbitrary point $\bar{x} \in  C$, parameter $\delta$ was taken as
\begin{eqnarray}
\begin{aligned}
& \delta := \bar{\delta} \| \bar{x}\|,
\end{aligned}
\end{eqnarray}
where $0.02 \leq \bar{\delta} \leq 0.10$ was also chosen at random.

Next, we refer to the method of dealing with the nonsmooth term $g_{j}$ in the subproblem of the generalized conditional gradient for MCOPs proposed by Assun\c{c}\~{a}o in \cite{AFP23} to deal with the nonsmooth term in the subproblem of our algorithm NPQNA for MCOPs.
We define $A_{j} = [B_{j}; -B_{j}] \in \mathbb{R}^{2n \times n}$ and $b_{j} := \delta e \in \mathbb{R}^{2n}$, then
\eqref{g_j} and \eqref{Z_j} can be rewritten as
\begin{eqnarray}
\begin{aligned}
& \ \ \max_{z} \langle x, z \rangle \\
& s.t. \ \ A_{j}z \preceq b_{j}.
\end{aligned}
\end{eqnarray}
Since the variables in (\ref{Optim-Model1}) are $t$ and $d$, we convert the nonsmooth term $g_{j}(x + d)$ into a dual form.
\begin{eqnarray}\label{g_j(x+d)}
\begin{aligned}
& \ \ \min_{w} \langle b_{j}, w \rangle \\
& s.t. \ \ A_{j}^{T}w = x + d, \\
& \ \ \ \ \ w \succeq 0.
\end{aligned}
\end{eqnarray}
By using duality theory, substituted the expression of (\ref{g_j(x+d)}) into (\ref{Optim-Model}), we can get
\begin{eqnarray}\label{dual}
\begin{aligned}
& \min \ t \\
& s.t. \ \nabla f_{j}(x)^{T} d + \frac{1}{2}d^{T}B_{j}(x)d + b_{j}^{T}w - g_{j}(x) - t \leq 0,\\
& \ \ \ \ \ \  A_{j}^{T}w_{j} = x + d, \\
& \ \ \ \ \ \  w_{j} \succeq 0, \  j \in \Lambda_{m}, \\
& \ \ \ \ \ \ (t,d) \in \mathbb{R} \times \mathbb{R}^{n}.
\end{aligned}
\end{eqnarray}
We briefly review the subproblem of the proximal Newton-type algorithm NPGA for MCOPs.
$$\min_{d \in \mathbb{R}^{n}}\max_{j \in \Lambda_{m}}\nabla f_{j}(x)^{T} d + \frac{1}{2}d^{T}\nabla^{2} f_{j}(x)d + g_{j}(x + d) - g_{j}(x).$$
Likewise, the subproblem of the NPGA for MCOPs can be reformulated as the
following quadratic programming problem
\begin{eqnarray}\label{MONPG_sub}
\begin{aligned}
& \min \ t \\
& s.t. \ \nabla f_{j}(x)^{T} d + \frac{1}{2}d^{T}\nabla^{2} f_{j}(x)d + b_{j}^{T}w - g_{j}(x) - t \leq 0,\\
& \ \ \ \ \ \  A_{j}^{j}w_{j} = x + d, \\
& \ \ \ \ \ \  w_{j} \succeq 0, \  j \in \Lambda_{m}, \\
& \ \ \ \ \ \ (t,d) \in \mathbb{R} \times \mathbb{R}^{n}.
\end{aligned}
\end{eqnarray}
And the subproblem of the proximal quasi-Newton algorithm PQNA for MCOPs.
$$\min_{d \in \mathbb{R}^{n}}\max_{j \in \Lambda_{m}}\nabla f_{j}(x)^{T} d + \frac{1}{2}d^{T}B_{j}(x)d + g_{j}(x + d) - g_{j}(x) + \frac{w}{2}\|d\|^{2}.$$
Likewise, the subproblem of the PQNA for MCOPs can be reformulated as the
following quadratic programming problem
\begin{eqnarray}\label{PQNA_sub}
\begin{aligned}
& \min \ t \\
& s.t. \ \nabla f_{j}(x)^{T} d + \frac{1}{2}d^{T}B_{j}(x)d  + \frac{w}{2}\|d\|^{2} + b_{j}^{T}w - g_{j}(x) - t \leq 0,\\
& \ \ \ \ \ \  A_{j}^{j}w_{j} = x + d, \\
& \ \ \ \ \ \  w_{j} \succeq 0, \  j \in \Lambda_{m}, \\
& \ \ \ \ \ \ (t,d) \in \mathbb{R} \times \mathbb{R}^{n}.
\end{aligned}
\end{eqnarray}
In our codes, we use the solver 'SLSQP' for solving convex quadratic programming problem in python 3.9 to solve problem (\ref{dual}), (\ref{MONPG_sub}) and (\ref{PQNA_sub}).

Pareto fronts using NPQNA, PQNA and NPGA of a two objective and a three
objective test problem (Problem 1, 3, 7, 9 and 14 in Table 1) are provided in Figures 1, 2, 3, 4 and 5
respectively. Similarly we have generated approximate Pareto fronts of all test problems
mentioned in Table 1.
\vspace{0.5em}
\begin{table}[H]
    \centering
    \caption{Details of test problems.}
    \label{tab:univ-compa}
    \resizebox{\textwidth}{45mm}{
    \begin{tabular}{llllllccccccrrrrrr}
    \hline
        \textbf{S1.No} & & & \textbf{(m,n)} & & & \textbf{f} & & & \textbf{lb$^{T}$} & & & \textbf{ub$^{T}$} & & & \textbf{Ref}  \\ \hline
        1 & & & $(2,2)$ & & & AN1 & & & (-3,-3) & & & (7,7) & & & \cite{AP15} \\
        2 & & & $(2,2)$ & & & AP2 & & & (-5,-5) & & & (5,5) & & & \cite{DAPM00} \\
        3 & & & $(2,2)$ & & & BK1 & & & (-3,-3) & & & (5,5) & & & \cite{CLV07} \\
        4 & & & $(3,3)$ & & & FDS & & & (-2,-2,-2) & & & (4,4,4) & & & \cite{CF98} \\
        5 & & & $(3,5)$ & & & FDS & & & (-2,...,-2) & & & (2,...,2) & & & \cite{CF98} \\
        6 & & & $(3,2)$ & & & IKK1 & & & (-2,-2) & & & (3,3) & & & \cite{CLV07} \\
        7 & & & $(3,3)$ & & & IKK1 & & & (-2,-2,-2) & & & (2,2,2) & & & \cite{CLV07} \\
        8 & & & $(2,2)$ & & & JOS1 & & & (-5,-5) & & & (5,5) & & & \cite{D17} \\
        9 & & & $(2,2)$ & & & LOVISON1 & & & (-3,-3) & & & (5,5) & & & \cite{DJ13} \\
        10 & & & $(2,2)$ & & & LRS1 & & & (-50,-50) & & & (50,50) & & & \cite{CLV07} \\
        11 & & & $(3,3)$ & & & MHHM2 & & & (-4,-4,-4) & & & (4,4,4) & & & \cite{CLV07} \\
        12 & & & $(3,2)$ & & & MHHM2 & & & (-4,-4) & & & (4,4) & & & \cite{CLV07} \\
        13 & & & $(2,1)$ & & & MOP1 & & & (-100) & & & (100) & & & \cite{CLV07} \\
        14 & & & $(3,3)$ & & & MOP7 & & & (-4,-4,-4) & & & (4,4,4) & & & \cite{CLV07} \\
        15 & & & $(3,2)$ & & & MOP7 & & & (-4,-4) & & & (4,4) & & & \cite{CLV07} \\
        16 & & & $(2,2)$ & & & MS1 & & & (-2,-2) & & & (2,2) & & & \cite{CCMP95} \\
        17 & & & $(2,4)$ & & & MS2 & & & (-2,...,-2) & & & (2,...,2) & & & \cite{CCMP95} \\
        18 & & & $(3,3)$ & & & SDD1 & & & (-2,-2,-2) & & & (2,2,2) & & & \cite{DGW92} \\
        19 & & & $(2,2)$ & & & SP1 & & & (-1,-1) & & & (5,5) & & & \cite{CLV07} \\
        20 & & & $(3,2)$ & & & VFM1 & & & (-2,-2) & & & (4,4) & & & \cite{CLV07} \\
        21 & & & $(3,3)$ & & & VFM1 & & & (-2,-2,-2) & & & (4,4,4) & & & \cite{CLV07} \\
        22 & & & $(2,2)$ & & & VU1 & & & (-3,-3) & & & (3,3) & & & \cite{CLV07} \\
        %23 & & & $(2,2)$ & & & VU2 & & & (-3,-3) & & & (3,3) & & & \cite{CLV07} \\
        23 & & & $(3,3)$ & & & ZLT1 & & & (-3,-3,-3) & & & (3,3,3) & & & \cite{CLV07} \\
    \hline
    \end{tabular}}
\end{table}
\vspace{0.5em}

\vspace{-0.8cm}
\begin{figure}[H]
		\centering
		\subfigure[Approximal Pareto fronts by NPQNA]
		{
			\begin{minipage}[H]{0.22\linewidth}
				\centering
                \includegraphics[scale=0.22]{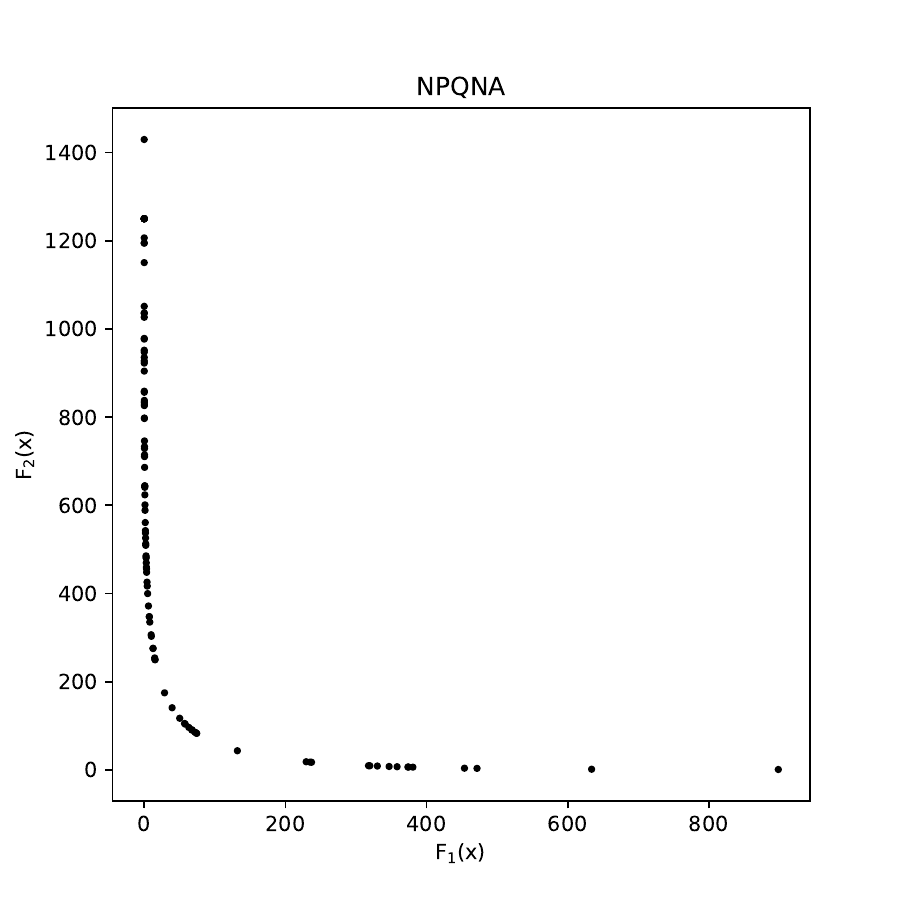} \\
			\end{minipage}
		}
\quad
		\subfigure[Approximal Pareto fronts by PQNA]
		{
			\begin{minipage}[H]{0.22\linewidth}
				\centering
				\includegraphics[scale=0.22]{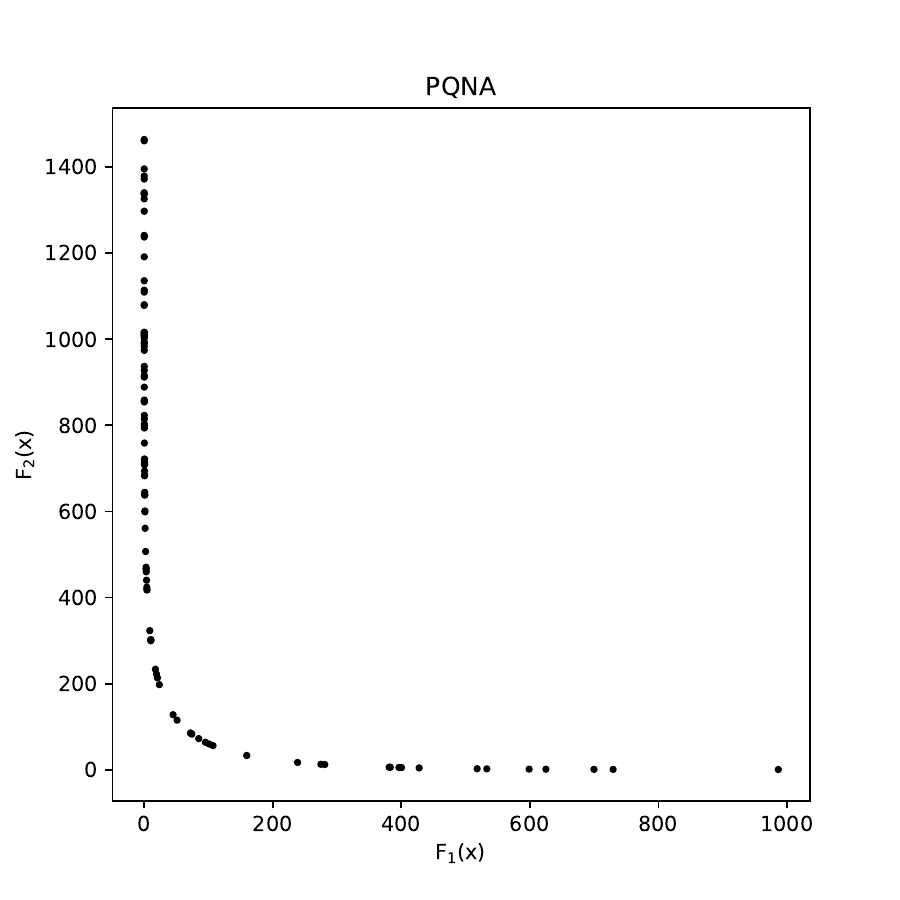} \\
			\end{minipage}
		}
\quad
		\subfigure[Approximal Pareto fronts by NPGA]
		{
			\begin{minipage}[H]{0.22\linewidth}
				\centering
				\includegraphics[scale=0.22]{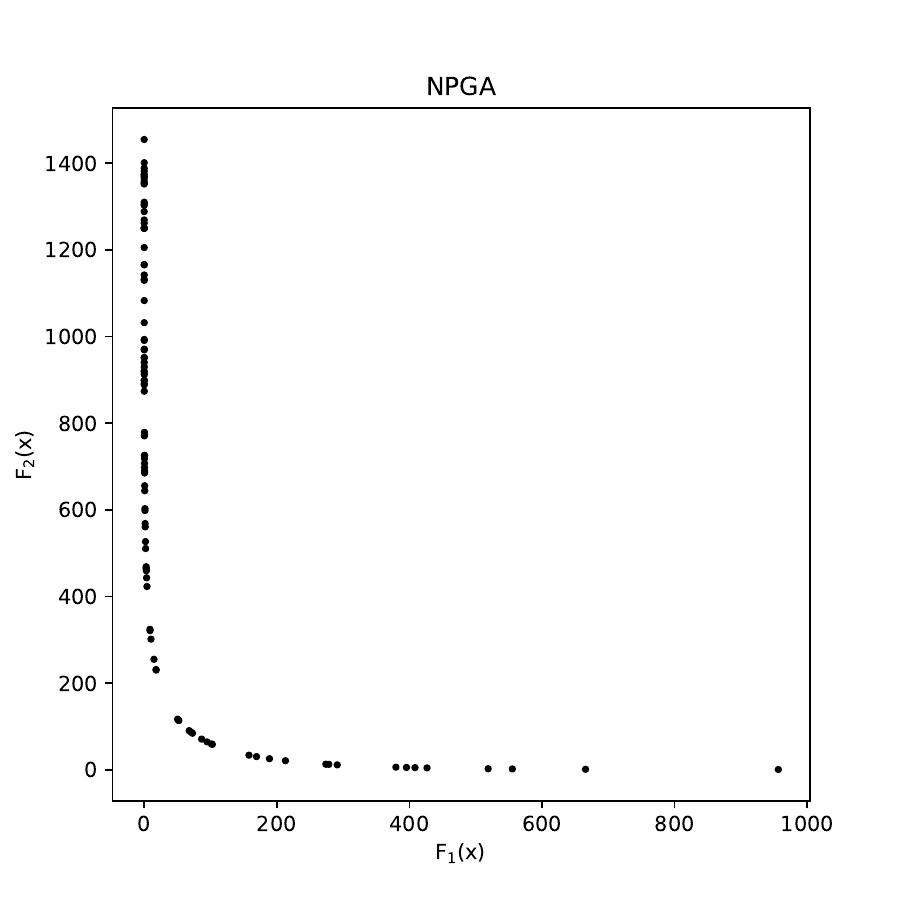} \\			
			\end{minipage}
		}
		\caption{Approximate Pareto fronts of problem 1. (a)Approximate Pareto fronts by NPQNA; (b)Approximate Pareto fronts by PQNA; (c)Approximate Pareto fronts by NPGA.}
		\label{f1}
\end{figure}
\vspace{-0.8cm}

\vspace{-0.8cm}
\begin{figure}[H]
		\centering
		\subfigure[Approximal Pareto fronts by NPQNA]
		{
			\begin{minipage}[H]{0.22\linewidth}
				\centering
                \includegraphics[scale=0.22]{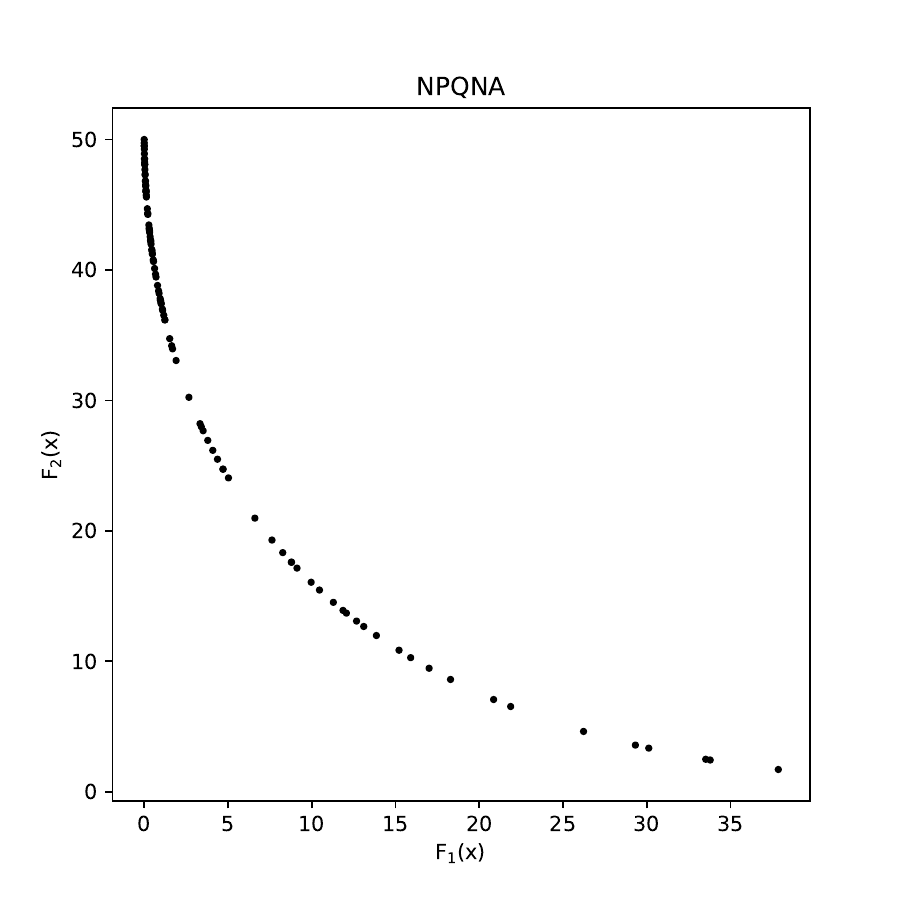} \\
			\end{minipage}
		}
\quad
		\subfigure[Approximal Pareto fronts by PQNA]
		{
			\begin{minipage}[H]{0.22\linewidth}
				\centering
				\includegraphics[scale=0.22]{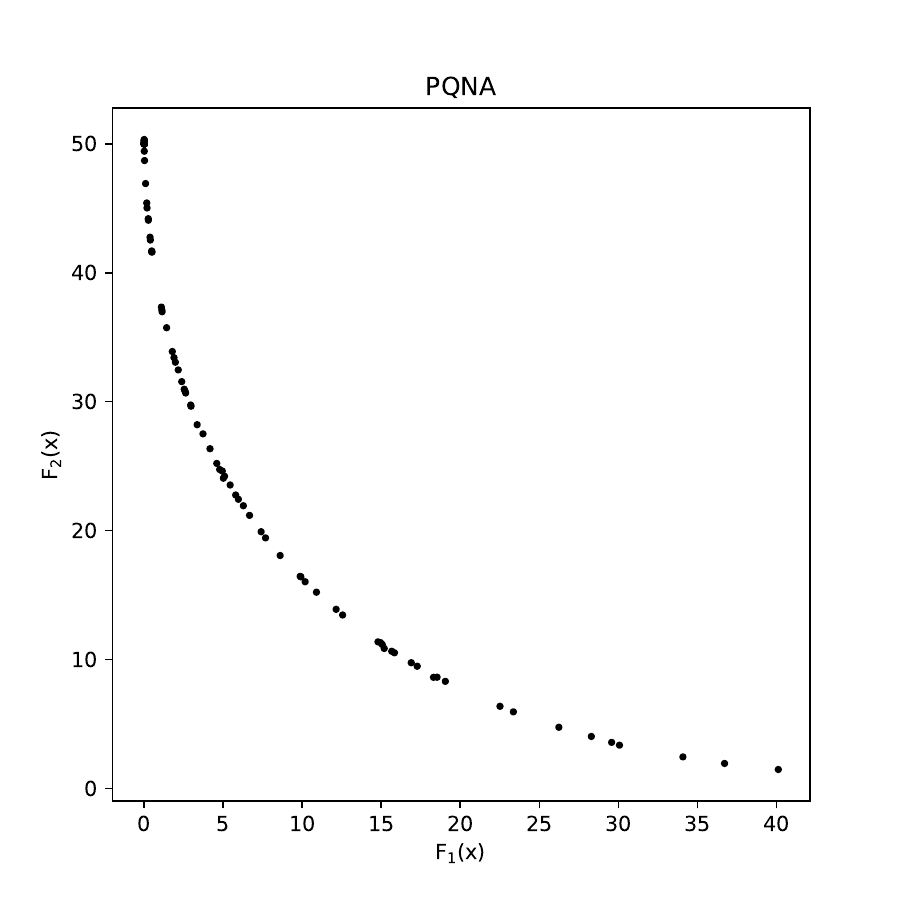} \\
			\end{minipage}
		}
\quad
		\subfigure[Approximal Pareto fronts by NPGA]
		{
			\begin{minipage}[H]{0.22\linewidth}
				\centering
				\includegraphics[scale=0.22]{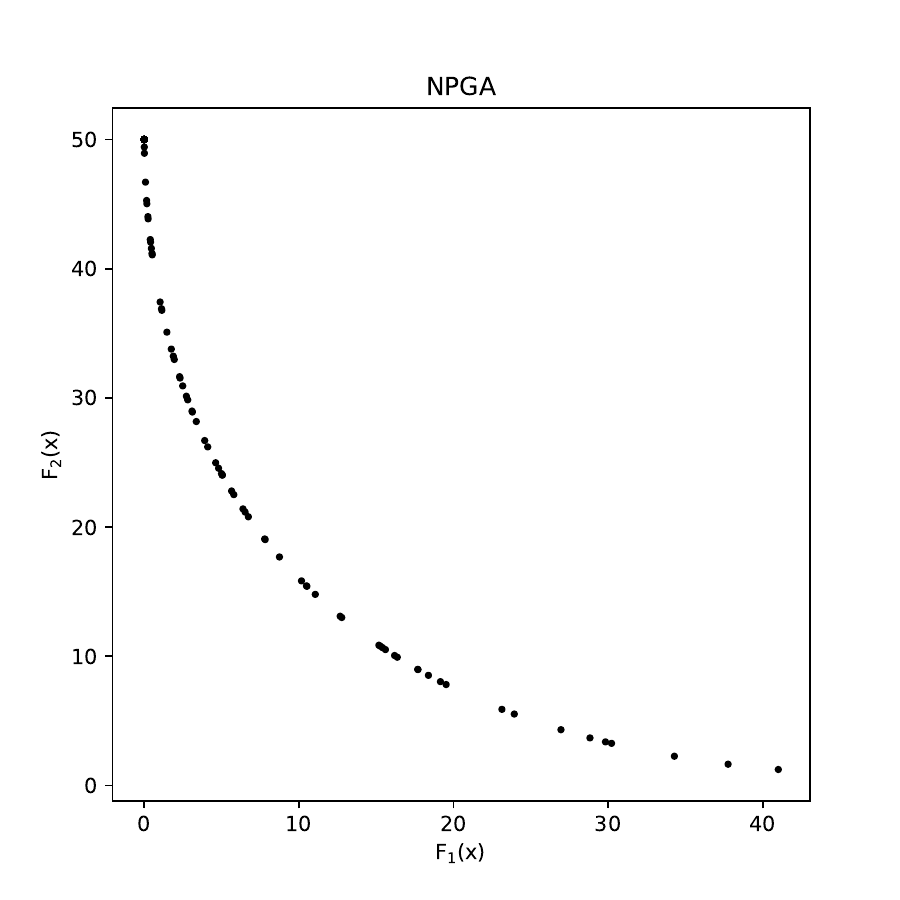} \\			
			\end{minipage}
		}
		\caption{Approximate Pareto fronts of problem 3. (a)Approximate Pareto fronts by NPQNA; (b)Approximate Pareto fronts by PQNA; (c)Approximate Pareto fronts by NPGA.}
		\label{f1}
\end{figure}
\vspace{-0.8cm}

\vspace{-0.8cm}
\begin{figure}[H]
		\centering
		\subfigure[Approximal Pareto fronts by NPQNA]
		{
			\begin{minipage}[H]{0.22\linewidth}
				\centering
                \includegraphics[scale=0.22]{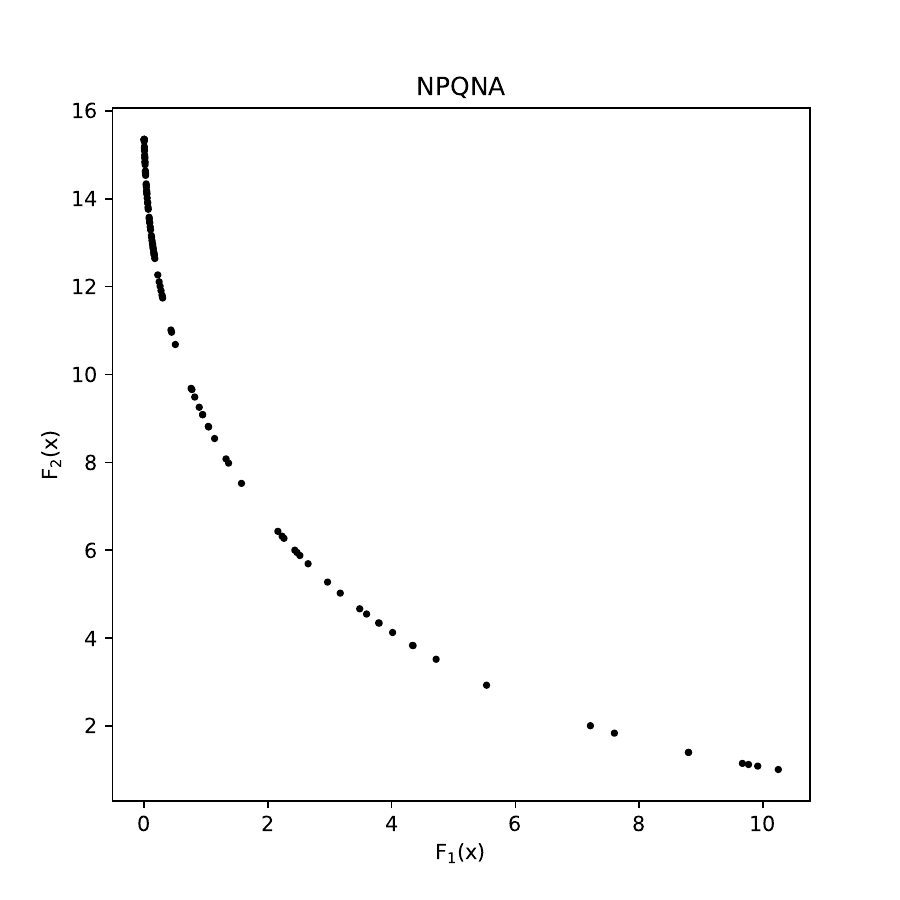} \\
			\end{minipage}
		}
\quad
		\subfigure[Approximal Pareto fronts by PQNA]
		{
			\begin{minipage}[H]{0.22\linewidth}
				\centering
				\includegraphics[scale=0.22]{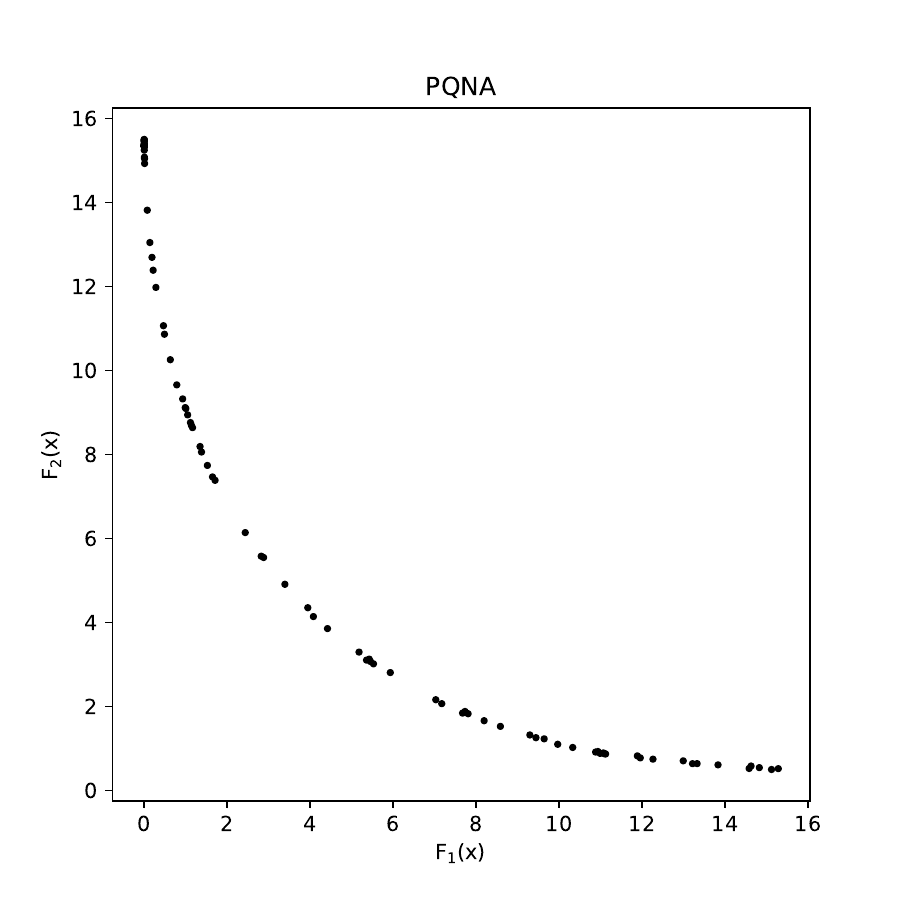} \\
			\end{minipage}
		}
\quad
		\subfigure[Approximal Pareto fronts by NPGA]
		{
			\begin{minipage}[H]{0.22\linewidth}
				\centering
				\includegraphics[scale=0.22]{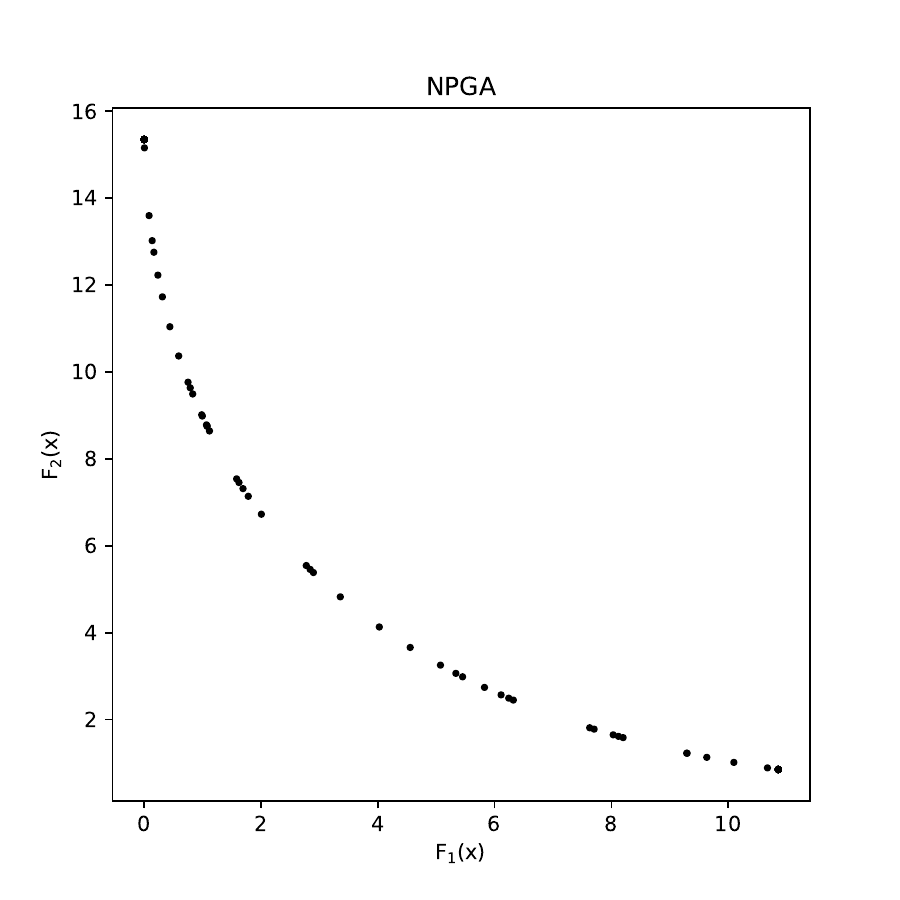} \\			
			\end{minipage}
		}
		\caption{Approximate Pareto fronts of problem 9. (a)Approximate Pareto fronts by NPQNA; (b)Approximate Pareto fronts by PQNA; (c)Approximate Pareto fronts by NPGA.}
		\label{f1}
\end{figure}
\vspace{-0.8cm}

\vspace{-0.8cm}
\begin{figure}[H]
		\centering
		\subfigure[Approximal Pareto fronts by NPQNA]
		{
			\begin{minipage}[H]{0.22\linewidth}
				\centering
                \includegraphics[scale=0.22]{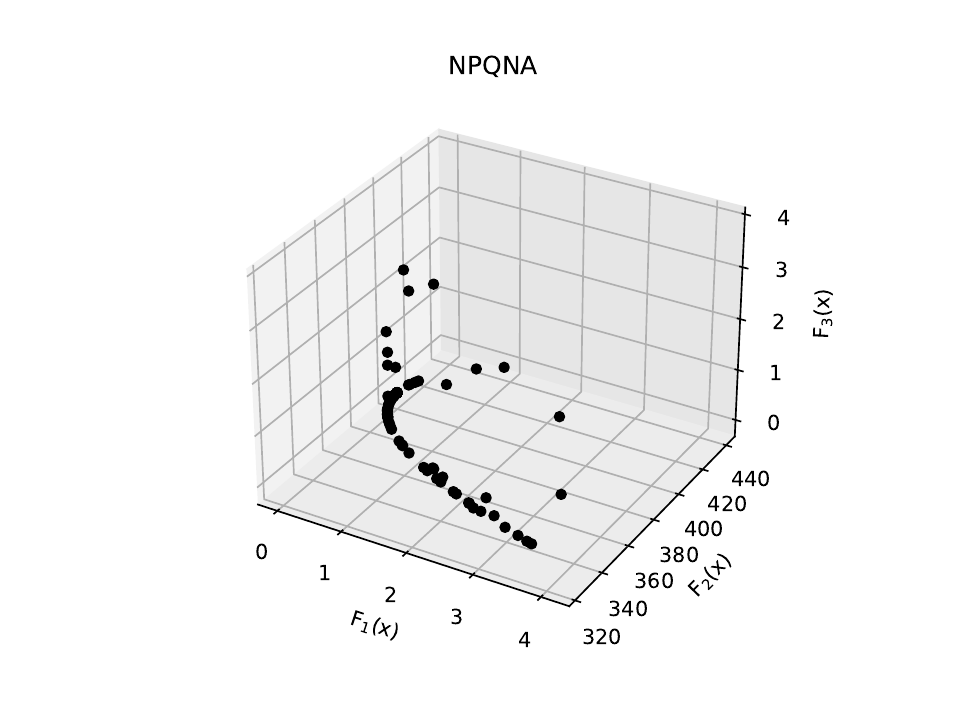} \\
			\end{minipage}
		}
\quad
		\subfigure[Approximal Pareto fronts by PQNA]
		{
			\begin{minipage}[H]{0.22\linewidth}
				\centering
				\includegraphics[scale=0.22]{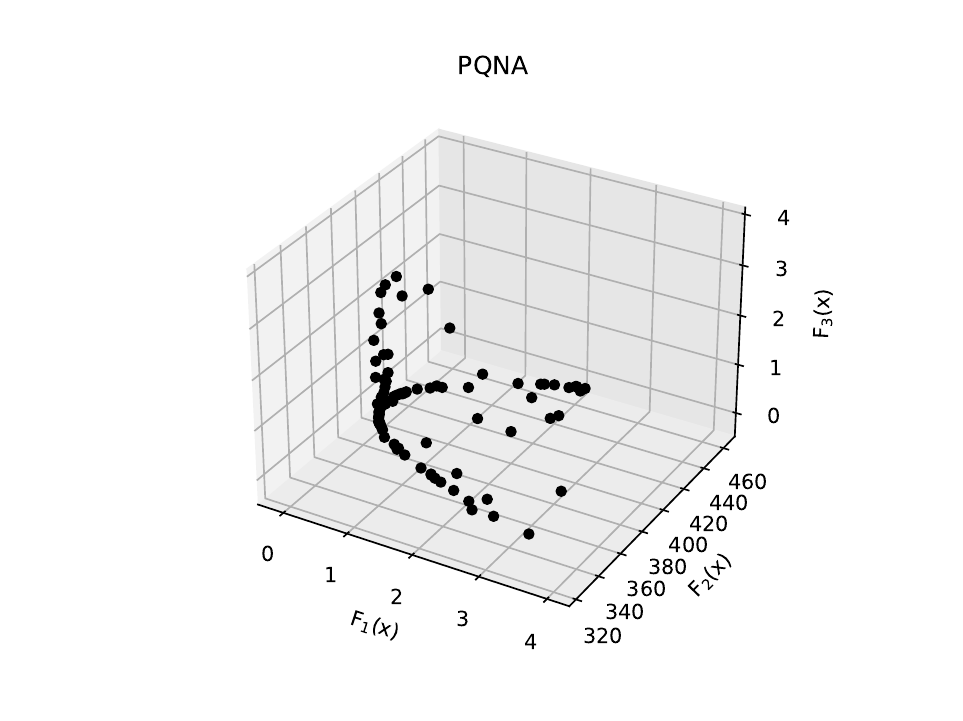} \\
			\end{minipage}
		}
\quad
		\subfigure[Approximal Pareto fronts by NPGA]
		{
			\begin{minipage}[H]{0.22\linewidth}
				\centering
				\includegraphics[scale=0.22]{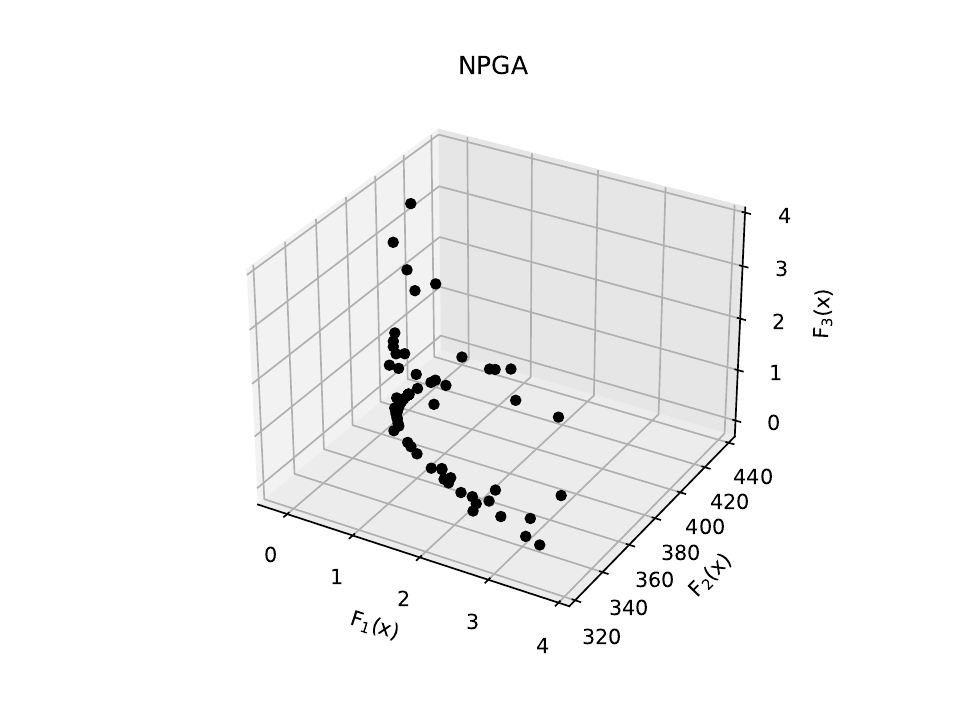} \\			
			\end{minipage}
		}
		\caption{Approximate Pareto fronts of problem 7. (a)Approximate Pareto fronts by NPQNA; (b)Approximate Pareto fronts by PQNA; (c)Approximate Pareto fronts by NPGA.}
		\label{f1}
\end{figure}
\vspace{-0.8cm}

\vspace{-0.8cm}
\begin{figure}[H]
		\centering
		\subfigure[Approximal Pareto fronts by NPQNA]
		{
			\begin{minipage}[H]{0.22\linewidth}
				\centering
                \includegraphics[scale=0.22]{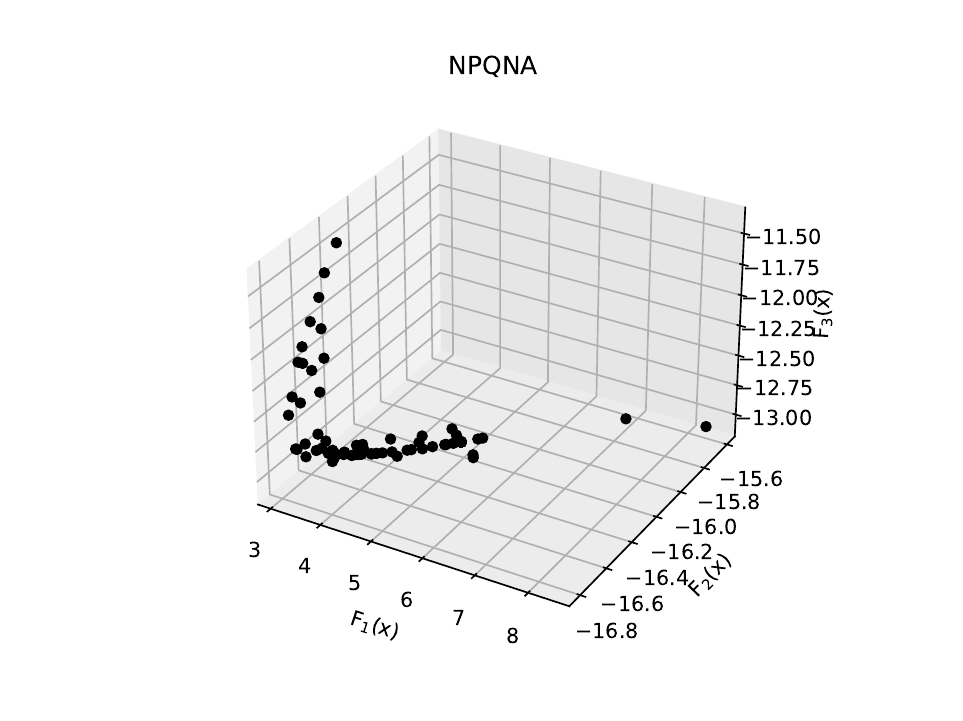} \\
			\end{minipage}
		}
\quad
		\subfigure[Approximal Pareto fronts by PQNA]
		{
			\begin{minipage}[H]{0.22\linewidth}
				\centering
				\includegraphics[scale=0.22]{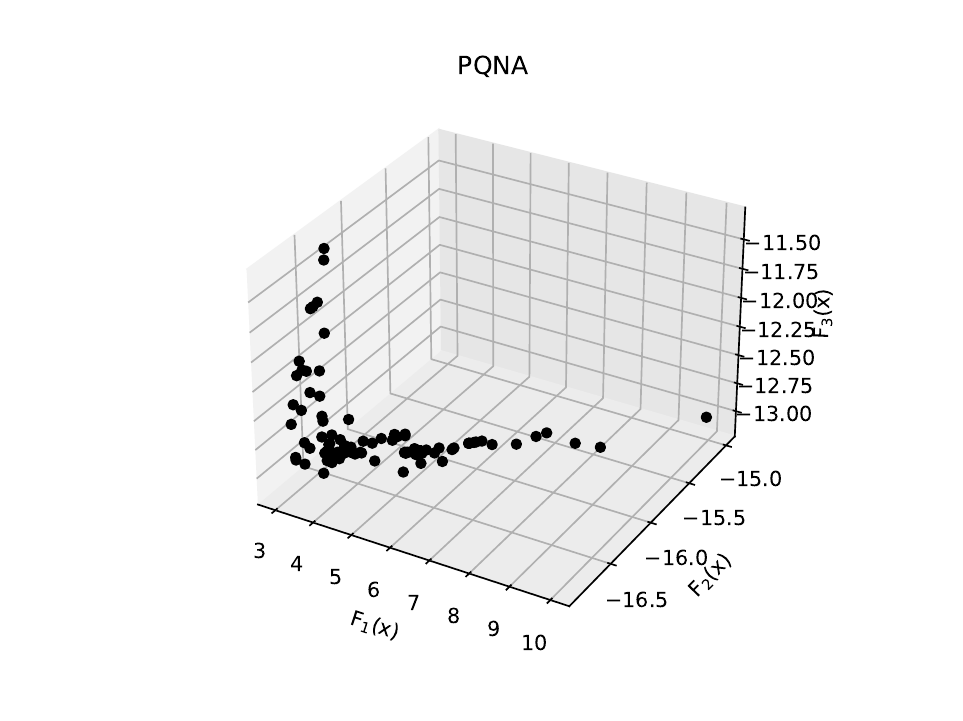} \\
			\end{minipage}
		}
\quad
		\subfigure[Approximal Pareto fronts by NPGA]
		{
			\begin{minipage}[H]{0.22\linewidth}
				\centering
				\includegraphics[scale=0.22]{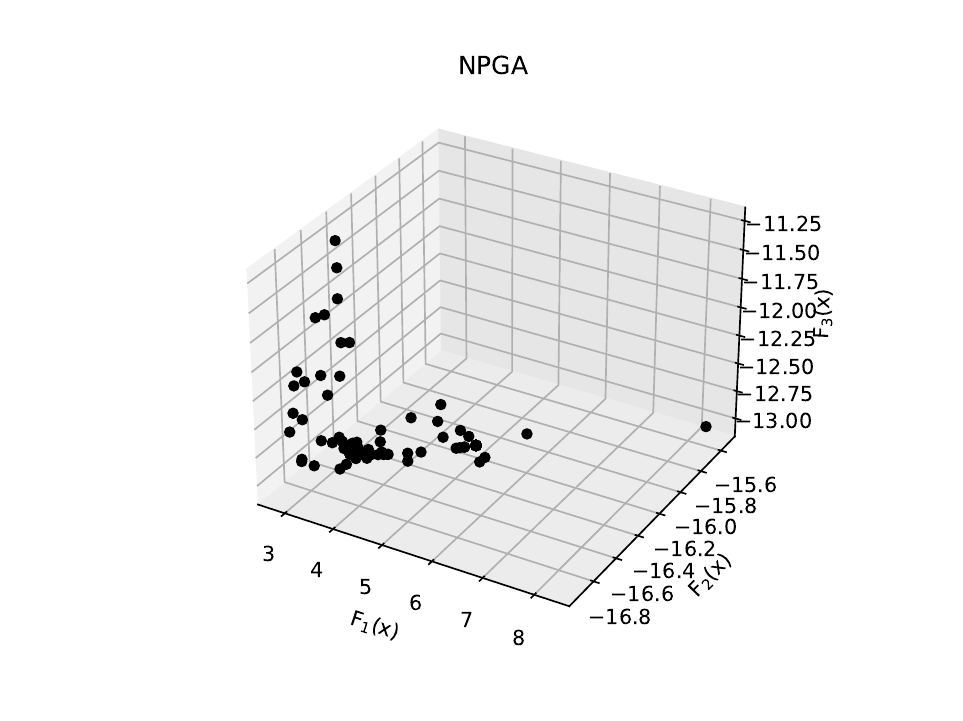} \\			
			\end{minipage}
		}
		\caption{Approximate Pareto fronts of problem 14. (a)Approximate Pareto fronts by NPQNA; (b)Approximate Pareto fronts by PQNA; (c)Approximate Pareto fronts by NPGA.}
		\label{f1}
\end{figure}
\vspace{-0.8cm}

\begin{table}[H]
    \centering
    \caption{Computational Details.}
    \label{tab:univ-compa}
    \resizebox{\textwidth}{43mm}{
    \begin{tabular}{*{7}{ccc}}
    \hline
    \multirow{2}*{Problem} & & \multicolumn{3}{c}{NPQNA} & & \multicolumn{3}{c}{PQNA} & & \multicolumn{3}{c}{NPGA} \\
    \cmidrule(llrr){2-5}\cmidrule(llrr){6-9}\cmidrule(llrr){10-13}
    & & it & f & CPU & & it & f & CPU & & it & f  & CPU  \\
    \hline
    1 & & \pmb{276.23} & \pmb{277.23} & \pmb{7.08} & & 300 & 301 & 13.28 & & 291.18 & 292.18 & 9.49 \\
    2 & & \pmb{161.88} & \pmb{162.88} & \pmb{7.95} & & 300 & 301 & 15.45 & & 300 & 301 & 11.32 \\
    3 & & 297.01 & 298.01 & 1.91 & & 226.34 & 227.34 & 3.23 & & \pmb{204.32} & \pmb{205.32} & \pmb{1.86} \\
    4 & & 177.53 & 178.53 & \pmb{1.34} & & 297.05 & 298.05 & 3.68  & & \pmb{162.52} & \pmb{163.52} & 1.35 \\
    5 & & 233.43 & 234.43 & 2.41 & & 297.07 & 298.07 & 4.5 & & \pmb{156.57} & \pmb{157.57} & \pmb{1.78} \\
    6 & & 249.3 & 250.3 & \pmb{1.96} & & 264.93 & 265.93 & 4.9 & & \pmb{225.25} & \pmb{226.25} & 2.71 \\
    7 & & 198.64 & 199.64 & \pmb{5.25} & & 270.78 & 271.78 & 12.73 & & \pmb{180.4} & \pmb{181.4} & 5.55 \\
    8 & & \pmb{198.34} & \pmb{199.34} & \pmb{0.79} & & 244.17 & 245.17 & 1.75 & & \pmb{198.34} & \pmb{199.34} & 0.96 \\
    9 & & 279.08 & 280.08 & \pmb{9.27} & & 235.02 & 236.02 & 16.99 & & \pmb{192.36} & \pmb{193.36} & 12.44 \\
    10 & & 255.24 & 256.24 & 2.78 & & \pmb{97.09} & \pmb{98.09} & 3.84 & & \pmb{96.68} & \pmb{97.68} & \pmb{2.33} \\
    11 & & \pmb{39.36} & \pmb{40.36} & \pmb{5.66} & & 252.39 & 253.39 & 15.74 & & 180.4 & 181.4 & 7.41 \\
    12 & & \pmb{50.47} & \pmb{51.47} & \pmb{2.39} & & 266.37 & 267.37 & 7.26 & & 162.97 & 163.97 & 3.92 \\
    13 & & \pmb{5.18} & \pmb{6.18} & \pmb{0.13} & & 165.25 & 166.25 & 3.32 & & 66.78 & 67.78 & 1.07 \\
    14 & & 201.59 & 202.59 & \pmb{4.05} & & 262.05 & 263.05 & 9.75 & & \pmb{180.4} & \pmb{181.4} & 4.50  \\
    15 & & \pmb{122.43} & \pmb{123.43} & \pmb{0.77} & & 294.18 & 295.18 & 2.82 & & 219.27 & 220.27 & 1.48 \\
    16 & & \pmb{146.28} & \pmb{147.28} & \pmb{5.6} & & 300 & 301 & 11.32 & & 228.24 & 229.24 & 7.8 \\
    17 & & \pmb{12.79} & \pmb{13.79} & \pmb{0.41} & & 300 & 301 & 10.43 & & 124.69 & 125.69 & 3 \\
    18 & & \pmb{239.84} & \pmb{240.84} & \pmb{9.38} & & 300 & 301 & 22.9 & & 297.08 & 298.08 & 13.04 \\
    19 & & \pmb{212.98} & \pmb{213.98} & \pmb{3.89} & & 297.79 & 298.79 & 6.03 & & 300 & 301 & 4.22 \\
    20 & & 214.78 & 215.78 & \pmb{3.48} & & 300 & 301 & 9.51 & & \pmb{102.66} & \pmb{103.66} & 4.6 \\
    21 & & \pmb{186.3} & \pmb{187.3} & \pmb{6.87} & & 300 & 301 & 18.98 & & 300 & 301 & 9.86 \\
    22 & & \pmb{300} & \pmb{301} & \pmb{4.75} & & \pmb{300} & \pmb{301} & 8.55 & & \pmb{300} & \pmb{301} & 5.39 \\
    %23 & & \pmb{1.7} & \pmb{2.7} & \pmb{4.77} & & 7.93 & 8.93 & 8.61 & & 300 & 301 & 6.37 \\
    23 & & \pmb{182.57} & \pmb{183.57} & \pmb{2.46} & & 300 & 301 & 6.34 & & 225.25 & 226.25 & 3.06 \\
    \hline
    \end{tabular}}
\end{table}

Numerical experiments show that our algorithm NPQNA of MCOPs has less iteration times, function evaluation times and CPU running time than NPGA of MCOPs in \cite{A22} and PQNA of MCOPs in \cite{PR22} on most test problems.
Although our algorithm NPQNA has more iterations than NPGA and PQNA on some test problems, our algorithm has less CPU running time, which also shows the superiority of our algorithm.
By analysing the related problems that the number of iterations of NPQNA is more than that of NPGA, we find that if the Hessian matrix of the objective function has an explicit expression, the number of iterations updated by Newton  is smaller than that updated by  quasi-Newton.
However, if the Hessian matrix of the objective function is not easy to obtain, our NPQNA algorithm has obvious advantages.

\section{Conclusions}

In this paper, we have developed  a nonmonotone proximal quasi-Newton algorithm for unconstrained convex MCOPs.
Under standard assumptions,  we proved that the sequence generated by this method converges to a Pareto optimal.
Furthermore, under the assumptions  that the Hessian
continuity, the Dennis-Mor\'{e} criterion and strong convexity, the local superlinear convergence rate of the algorithm is obtained.
We consider convex MCOPs, but in real life, many MCOPs are nonconvex, so it is of practical value to study the algorithm for solving nonconvex MCOPs.
And this method only starts from the quadratic approximation of the smooth function, and does not take full advantage of the characteristics of the nonsmooth term.
Lin et al.\cite{LMH19} proposed a proximal quasi-Newton algorithm in scalar optimization, which makes full use of the properties of nonsmooth function, performs Moreau envelope processing on nonsmooth term, and then uses the properties of Moreau envelope to perform quasi-Newton acceleration on the whole smooth problem.
This is a very nice approach, but there is no research on MCOPs yet.
In the future, we try to study smoothing techniques for MCOPs and nonconvex MCOPs.

%References

\end{document}